\documentclass[12pt,letterpaper,reqno]{amsart}
\usepackage{mathrsfs}
\usepackage{times}
\usepackage[T1]{fontenc}
\usepackage[reqno]{amsmath}
\usepackage{amsfonts,amsthm,amssymb,amscd}

\renewcommand{\contentsname}{Contents\\{\footnotesize\normalfont(A table
of contents should normally not be included)}}

\newcommand\be{\begin{equation}}
\newcommand\ee{\end{equation}}
\newcommand\bea{\begin{eqnarray}}
\newcommand\eea{\end{eqnarray}}
\newcommand\bi{\begin{itemize}}
\newcommand\ei{\end{itemize}}
\newcommand\ben{\begin{enumerate}}
\newcommand\een{\end{enumerate}}
\newcommand\bc{\begin{center}}
\newcommand\ec{\end{center}}
\newcommand\ba{\begin{array}}
\newcommand\ea{\end{array}}

\newtheorem{thm}{Theorem}[section]

\newtheorem{cor}[thm]{Corollary}
\newtheorem{lem}[thm]{Lemma}

\newtheorem{defi}[thm]{Definition}

\def\notdiv{\ \mathbin{\mkern-8mu|\!\!\!\smallsetminus}}
\newcommand{\Qoft}{\Bbb{Q}(t)}  
\newcommand{\done}{\Box} 
\newcommand{\R}{\ensuremath{\Bbb{R}}}
\newcommand{\C}{\ensuremath{\Bbb{C}}}
\newcommand{\Z}{\ensuremath{\Bbb{Z}}}
\newcommand{\Q}{\Bbb{Q}}

\newcommand{\foh}{\frac{1}{2}}  

\newcommand{\fologn}{O\Big(\frac{1}{\log N}\Big)}

\newcommand{\flogpn}{\frac{\log p}{\log N}}

\newcommand{\js}[1]{ { \underline{#1} \choose p} }

\newcommand{\oof}{\frac{1}{|\mathcal{F}|}} 
\newcommand{\oop}{\frac{1}{p}}
\newcommand{\oopi}{\frac{1}{p_i}}

\newcommand{\hfi}{\widehat{f_i}}
\newcommand{\hfo}{\widehat{f_1}}
\newcommand{\hft}{\widehat{f_2}}
\newcommand{\NF}{|\mathcal{F}|}   
\newcommand{\sumtinF}{\sum_{t \in \mathcal{F}}}
\newcommand{\ooNF}{\frac{1}{|\mathcal{F}|}}   
\newcommand{\sumef}{\sum_{E \in \mathcal{F}}}

\newcommand{\plogx}{ \frac{\log p}{\log X}}
\newcommand{\plogm}{\frac{\log p}{\log M}}

\newcommand{\plogne}{\frac{\log p}{\log N_E}}

\newcommand{\pilogne}{\frac{\log p_i}{\log N_E}}

\newcommand{\kkot}[1]{ \frac{\sin \pi {#1} }{\pi {#1} } }

\begin{document}

\title[$1$- and $2$-Level Densities for Rational Families of Elliptic
Curves]{$1$- and $2$-Level Densities for Rational Families of
Elliptic Curves: Evidence for the Underlying Group Symmetries}
\author{Steven J. Miller}
\email{sjmiller@math.princeton.edu}
\address{Department of Mathematics, Princeton University, Princeton, NJ $08544$, U.S.A.}

\begin{abstract}
Following Katz-Sarnak \cite{KS1}, \cite{KS2}, Iwaniec-Luo-Sarnak
\cite{ILS}, and Rubinstein \cite{Ru}, we use the $1$- and
$2$-level densities to study the distribution of low lying zeros
for one-parameter rational families of elliptic curves over
$\Bbb{Q}(t)$. Modulo standard conjectures, for small support the
densities agree with Katz and Sarnak's predictions. Further, the
densities confirm that the curves' $L$-functions behave in a
manner consistent with having $r$ zeros at the critical point, as
predicted by the Birch and Swinnerton-Dyer conjecture. By studying
the $2$-level densities of some constant sign families, we find
the first examples of families of elliptic curves where we can
distinguish $SO(\mbox{even})$ from $SO(\mbox{odd})$ symmetry.
\end{abstract}

\maketitle



\section{Introduction}
\setcounter{equation}{0}

\subsection{$n$-Level Correlations and Densities}

Assuming GRH, the zeros of any $L$-function lie on the critical
line, and therefore it is possible to investigate statistics of
the normalized zeros. The general philosophy, born out in many
examples (see \cite{CFKRS}), is that the behavior of random
matrices / ensembles of random matrices behave similar to that of
$L$-functions / families of $L$-functions. By a family
$\mathcal{F}$ we mean a collection of geometric objects and their
associated $L$-functions, where the geometric objects have similar
properties.

We expect there is a symmetry group $\mathcal{G}(\mathcal{F})$
(one of the classical compact groups $U(N)$, $SU(N)$, $USp(2N)$,
$SO(\mbox{even})$ and $SO(\mbox{odd})$) which can be associated to
a family of $L$-functions, and that the behavior of eigenvalues of
matrices in $\mathcal{G}(\mathcal{F})$ should (after appropriate
normalizations) equal the behavior of zeros of $L$-functions.

Iwaniec, Luo and Sarnak \cite{ILS} consider (among other examples)
all cuspidal newforms of a given level and weight. Rubinstein
\cite{Ru} considers twists by fundamental discriminants $D$ of a
fixed modular form.

We study the family of all elliptic curves and various
one-parameter families of elliptic curves. Thus, in our case the
notion of family is the standard one from geometry: we have a
collection of curves over a base, and the geometry is much clearer
in our examples than in \cite{ILS} and \cite{Ru}.

Let $\{\alpha_j\}$ be an increasing sequence of numbers tending to
infinity, such as eigenvalues or zeros normalized to have mean
spacing $1$. For a compact box $B \subset \R^{n-1}$, define the
\emph{$n$-level correlation} by

\begin{eqnarray}
\lim_{N \to \infty} \frac{\# \Big\{ (\alpha_{j_1}-\alpha_{j_2},
\dots, \alpha_{j_{n-1}} - \alpha_{j_n}) \in B, j_i \in
\{1,\dots,N\}, j_i \neq j_k \Big\}}{N} \end{eqnarray}

Note that the $n$-level correlations are unaffected by removing
finitely many zeros. Instead of using a box, one can study a
smoothed version with a test function on $\R^n$ (see \cite{RS}).

For test functions whose Fourier Transform has small support,
Montgomery \cite{Mon} proved the $2$- and Hejhal \cite{Hej} proved
the $3$-level correlations for the zeros of $\zeta(s)$ are the
same as that of the GUE, and Rudnick-Sarnak \cite{RS} proved the
$n$-level correlations for all automorphic cuspidal $L$-functions
are the same as the GUE. The universality is due to the fact that
the correlations are controlled by the second moment of the
$a_p$'s, and while there are many possible limiting distributions,
all have the same second moment.

Katz and Sarnak \cite{KS1} prove the classical compact groups have
the same $n$-level correlations. In particular, we cannot use the
$n$-level correlations to distinguish GUE behavior, $U(N)$, from
the other classical compact groups. We are led to investigate
another statistic which will depend on the underlying group.

For $L$-functions of elliptic curves, the order of vanishing of
$L(s,E)$ at $s = \foh$ is conjecturally equal to the geometric
rank of the Mordell-Weil group (Birch and Swinnerton-Dyer
conjecture). If we force the Mordell-Weil group to be large, we
expect many zeros exactly at $s = \foh$, and this might influence
the behavior of the neighboring zeros. Hence we are led to study
the distribution of the first few, or low lying, zeros, and the
fascinating possibility that there could be a difference in
statistics for zeros near $\foh$ than for zeros higher up.

Let $f(x)$ be an even Schwartz function whose Fourier Transform is
supported in a neighborhood of the origin. We assume $f$ is of the
form $\prod_{i=1}^n f_i(x_i)$. The \emph{$n$-level density} for
the family $\mathcal{F}$ with test function $f$ is

\begin{eqnarray}
D_{n,\mathcal{F}}(f) = \oof \sumef \sum_{j_1,\dots, j_n \atop j_i
\neq \pm j_k} f_1\Big(\frac{\log
N_E}{2\pi}\gamma_E^{(j_1)}\Big)\cdots f_n\Big(\frac{\log
N_E}{2\pi}\gamma_E^{(j_n)}\Big), \end{eqnarray}

where $\gamma_E^{(j_i)}$ runs through the non-trivial zeros of the
curve $E$, and $N_E$ is its conductor. We rescale the zeros by
$\log N_E$ as this is the order of the number of zeros with
imaginary part less than a large absolute constant (see
\cite{ILS}). As $f_i$ is Schwartz, most of the contribution is due
to the zeros near the critical point. We use the Explicit Formula
(Equation \ref{thmef}) to relate sums of test functions over zeros
to sums over primes of $a_E(p)$ and $a_E^2(p)$.

Katz and Sarnak \cite{KS1} determine the $N \to \infty$ limits for
the $n$-level densities of eigenvalues near $1$ for the classical
compact groups (see Section \ref{secksdensities}); their
calculations can be modified to determine the densities of
classical compact groups with a forced number of eigenvalues at
$1$. Forcing eigenvalues at $1$ corresponds to $L$-functions with
zeros forced at the critical point.

\subsection{Results}

To any geometric family in the function field case, the results of
Katz and Sarnak (\cite{KS1}, \cite{KS2}) state the $n$-level
density of zeros near $\foh$ depends only on a symmetry group
attached to the family. In particular, for generic families of
elliptic curves the relevant symmetry is orthogonal. One can
further analyze the distributions depending on the signs of the
functional equations. As the families of elliptic curves are
self-dual, we expect the densities to be controlled by the
distribution of signs (all even: $SO(\mbox{even})$; all odd:
$SO(\mbox{odd})$; equidistributed: $O$).

For an elliptic curve $E_t$, let $D(t)$ be the product of the
irreducible polynomial factors of the discriminant $\Delta(t)$,
and let $C(t)$ be the conductor. Let $B$ be the largest square
dividing $D(t)$ for all $t$. Pass to a subsequence $ct + t_0$, and
call $t \in [N,2N]$ \emph{good} if $D(ct+t_0)$ is square-free,
except for primes $p|B$ where the power of such $p|D(t)$ is
independent of $t$.

The main result is Theorem \ref{thmonetwodensityrank}: \\

\textbf{Rational Surfaces Density Theorem:} \emph{Consider a
one-parameter family of elliptic curves of rank $r$ over $\Qoft$
which constitutes a rational surface. Assume GRH, $j(E_t)$
non-constant, and if $\Delta(t)$ has an irreducible polynomial
factor of degree at least $4$, assume the ABC Conjecture. }

\emph{After passing to a subsequence, for $t$ good, $C(t)$ is a
polynomial. Let $f_i$ be an even Schwartz function of small but
non-zero support $\sigma_i$ ($\sigma_1 < \min(\foh,\frac{2}{3m})$
for the $1$-level density, $\sigma_1 + \sigma_2 < \frac{1}{3m}$
for the $2$-level density).}

\emph{The $1$-level density agrees with the orthogonal densities
plus a term which equals the contributions from $r$ zeros at the
critical point. The $2$-level density agrees with
$SO(\mbox{even})$, $O$, and $SO(\mbox{odd})$ depending on whether
the signs are all even, equidistributed in the limit, or all odd,
plus a term which equals the contribution from $r$ zeros at the
critical point. Thus, for small support, the densities of the
zeros agree with Katz and Sarnak's predictions. Further, the
densities confirm that the curves' $L$-functions behave in a
manner consistent with having $r$ zeros at the critical point, as
predicted by the Birch and Swinnerton-Dyer conjecture.} \\

The ABC Conjecture is used to handle large prime divisors of
polynomials of degree $4$ or more (see \cite{Gr}). In place of
ABC, one could assume the Square-Free Sieve Conjecture.

For the $1$-level densities, the three orthogonal densities agree
for test functions with support less than $1$, split (ie, are
distinguishable) for support greater than $1$, but are all
distinguishable from $U$ and $Sp$ for any support. Hence, unlike
the $n$-level correlations, the $1$-level density is already
sufficient to observe non-GUE and non-symplectic behavior.

The polynomial growth of the conductor in families of elliptic
curves makes it difficult to evaluate the sums over primes for
test functions with moderate support. Converting to our language,
for small support the $1$-level densities for many families have
been shown equal to the Katz-Sarnak predictions: all elliptic
curves (Brumer and Heath-Brown \cite{Br}, \cite{BHB5}, support
less than $\frac{2}{3}$); twists of a given curve (support less
than $1$); one-parameter families (Silverman \cite{Si3}, small
support).

None of these are sufficient to distinguish the three orthogonal
candidates. Further, previous investigations have rescaled each
curve's zeros by the average of the logarithms of the conductors.
This greatly simplifies the calculations; however, the
normalization is no longer natural for each curve, as each curve
can sit in infinitely many families, each with a different average
spacing. By using local normalizations for each curve's zeros, the
$n$-level density for a family becomes the average of the
$n$-level densities for each curve.

The utility of the $2$-level density is that, even for test
functions with arbitrarily small support, the three candidate
orthogonal symmetries \emph{are} distinguishable, and in a very
satisfying way. The three candidates differ by a factor which
encodes the distribution of sign in the family, and all differ
from the GUE's $2$-level density.

We will study several families of constant sign, and we will see
that the densities are as expected. Thus, for these constant sign
families, the $2$-level density reflects the predicted symmetry,
which is invisible through the $1$-level density because of
support considerations.

Similar to the universality Rudnick and Sarnak \cite{RS} found in
studying $n$-level correlations, our universality follows from the
sums of $a_t^2(p)$ in our families (the second moments). For
non-constant $j(E_t)$, this follows from a Sato-Tate law proved by
Michel \cite{Mi} (Theorem \ref{thmmichel}).

\subsection{Structure of the Paper}

First, we calculate sums of the Fourier coefficients of elliptic
curves. We quote the predicted densities, and then calculate
useful expansions for the $1$- and $2$-level densities for
families of elliptic curves over $\Qoft$. We derive the density
results, conditional on the evaluation of many elliptic curve
sums. We calculate these sums for one-parameter rational families
of elliptic curves. We conclude with several examples (four
constant sign families, a rank $1$ and a rank $6$ rational
family).

We need excellent control over the conductors to evaluate the
above sums; the estimation is so delicate that if the log of
conductors are of size $m \log N$, fluctuations of size $O(1)$
yield error terms greater than the expected main terms.

The key observation is that the error terms can be controlled if
the conductors are monotone. By straightforward sieving and
applications of Tate's algorithm (to calculate the conductors),
given a one-parameter rational family of elliptic curves, we may
pass to a positive percent sub-family where the conductors are
monotone. Proofs of these results are given in the appendices.

In this paper, we concentrate on rational elliptic surfaces,
because here Tate's conjecture is known. Rosen and Silverman
\cite{RSi} show Tate's conjecture implies certain sums over primes
are related to the rank of the family over $\Qoft$. This will
allow us to interpret some of our density terms as the
contributions from $r$ critical point zeros.

The modifications needed to handle the family of all elliptic
curves, parametrized by

\be y^2 \ = \ x^3 + ax + b, \ \ a \in [-N^2,N^2], \ b \in
[-N^3,N^3], \ee

are straightforward, and can be found in \cite{Mil}.

Finally, if instead we normalize by the average of the logarithms
of the conductors, we obtain the same results, but with
significantly less work. This is done for one-parameter families
and the family of all elliptic curves in \cite{Mil}.


\section{Elliptic Curve Preliminaries}
\setcounter{equation}{0}

\subsection{Definitions}

Consider a one-parameter family $\mathcal{E}$ of elliptic curves
$E_t$ over $\Q(t)$:

\begin{eqnarray}
\mathcal{E}: y^2 + a_1(t)xy + a_3(t)y = x^3 + a_2(t)x^2 + a_4(t)x
+ a_6(t), \ \ a_i(t) \in \Z[t]. \end{eqnarray}

For each curve $E_t$, let $\Delta(t)$ be its discriminant and
$C(t)$ its conductor. Let $D(t)$ denote the product of the
irreducible polynomial factors dividing $\Delta(t)$. We will take
$t \in [N,2N]$ such that $D(t)$ is square-free.

Let $a_t(p) = a_{E_t}(p) = p+1 - N_{t,p}$, where $N_{t,p}$ is the
number of solutions of $E_t$ mod $p$ (including $\infty$). If $y^2
= x^3 + A(t)x + B(t)$, then

\be a_t(p) \ = \ -\sum_{t (p)} \js{x^3 + A(t)x + B(t)}. \ee

\subsection{Assumptions}

We assume the following at various points: \\

\textbf{Generalized Riemann Hypothesis (for Elliptic Curves)}
\emph{Let $L(s,E)$ be the (normalized) $L$-function of an elliptic
curve $E$. The non-trivial zeros $\rho$ of $L(s,E)$ have
$\mbox{Re}(\rho) = \foh$.} \\

Occasionally we assume the RH for the Riemann Zeta-function and Dirichlet
$L$-functions. \\

\textbf{Birch and Swinnerton-Dyer Conjecture \cite{BSD1},
\cite{BSD2}} \emph{Let $E$ be an elliptic curve of geometric rank
$r$ over $\Q$ (the Mordell-Weil group is $\Z^r \oplus T$). Then
the analytic rank (the order of vanishing of the $L$-function at
the critical point) is also $r$.} \\

We assume the above only for interpretation purposes. \\

\textbf{Tate's Conjecture for Elliptic Surfaces \cite{Ta}}
\emph{Let $\mathcal{E}/ \Q$ be an elliptic surface and
$L_2(\mathcal{E},s)$ be the $L$-series attached to
$H^2_{\mbox{{\'e}t}}(\mathcal{E}/ \overline{\Q}, \Q_l)$.
$L_2(\mathcal{E},s)$ has a meromorphic continuation to $\C$ and
$-\mbox{ord}_{s=1} L_2(\mathcal{E},s)$ $= \mbox{rank}\
NS(\mathcal{E}/ \Q)$, where $NS(\mathcal{E}/ \Q)$ is the
$\Q$-rational part of the N{\'e}ron-Severi group of $\mathcal{E}$.
Further, $L_2(\mathcal{E},s)$ does not vanish on
the line $\mbox{Re}(s) = 1$.} \\

Most of the one-parameter families that we investigate are
rational surfaces, in which case Tate's conjecture is known (see
 \cite{RSi}). \\

\textbf{ABC Conjecture} \emph{Fix $\epsilon > 0$. For co-prime
positive integers $a$, $b$ and $c$ with $c = a+b$ and $N(a,b,c) =
\prod_{p|abc} p$, $c \ll_\epsilon N(a,b,c)^{1+\epsilon}$.} \\

The full strength of ABC is never needed; rather, we need a
consequence of ABC, the Square-Free Sieve (see \cite{Gr}): \\

\textbf{Square-Free Sieve Conjecture} \emph{Fix an irreducible
polynomial $f(t)$ of degree at least $4$. As $N \to \infty$, the
number of $t \in [N,2N]$ with $f(t)$ divisible by $p^2$ for some
$p > \log N$ is $o(N)$.} \\

For irreducible polynomials of degree at most $3$, the above is
known, complete with a better error than $o(N)$ (\cite{Ho},
chapter $4$).

We use the Square-Free Sieve to handle the variations in the
conductors. If our evaluation of the log of the conductors is off
by as little as a small constant, the prime sums become
untractable. This is why many works normalize by the average
log-conductor. \\

\textbf{Restricted Sign Conjecture (for the Family $\mathcal{F}$)}
\emph{Consider a one-parameter family $\mathcal{F}$ of elliptic
curves. As $N \to \infty$, the signs of the curves $E_t$ are
equidistributed for $t \in [N,2N]$.} \\

The Restricted Sign conjecture often fails. First, there are
families with constant $j(E_t)$ where all curves have the same
sign.

Helfgott \cite{He} has recently related the Restricted Sign
conjecture to the Square-Free Sieve conjecture and standard
conjectures on sums of Moebius: \\

\textbf{Polynomial Moebius} \emph{Let $f(t)$ be a non-constant
polynomial such that no fixed square divides $f(t)$ for all $t$.
Then $\sum_{t=N}^{2N} \mu(f(t)) = o(N)$.} \\

The Polynomial Moebius conjecture is known for linear $f(t)$.

Helfgott shows the Square-Free Sieve and Polynomial Moebius imply
the Restricted Sign conjecture for many families. More precisely,
let $M(t)$ be the product of the irreducible polynomials dividing
$\Delta(t)$ and not $c_4(t)$. \\

\textbf{Theorem: Equidistribution of Sign in a Family \cite{He}:}
\emph{Let $\mathcal{F}$ be a one-parameter family with $a_i(t) \in
\Z[t]$. If $j(E_t)$ and $M(t)$ are non-constant, then the signs of
$E_t$, $t \in [N,2N]$, are equidistributed as $N \to \infty$.
Further, if we restrict to good $t$, $t \in [N,2N]$ such that
$D(t)$ is good (usually square-free), the signs are still
equidistributed in the limit.}

The above is only used to calculate $N(\mathcal{F},-1)$, the
percent of odd curves. Without this, we can still calculate the
$1$-level densities for small support, and all but one term in the
$2$-level densities, $N(\mathcal{F},-1) f_1(0)f_2(0)$.

\subsection{Explicit Formula}

The starting point for working with zeroes of the $L$-functions of
elliptic curves is the Explicit Formula (see \cite{Mes}), which
relates sums over zeros to sums over primes.

For an elliptic curve $E$ with conductor $N_E$,

\begin{eqnarray}\label{thmef}
\sum_{\gamma_E^{(j)}} G\Big(\gamma_E^{(j)} \frac{\log
N_E}{2\pi}\Big) & = & \widehat{G}(0) + G(0) -  2 \sum_p \plogne
\oop \widehat{G} \Big( \plogne \Big) a_E(p) \nonumber\\ & & - 2
\sum_p \plogne \frac{1}{p^2} \widehat{G} \Big(\frac{2 \log p}{\log
N_E} \Big) a_E^2(p) \nonumber\\ & & + O\Big(\frac{\log \log
N_E}{\log N_E}\Big).
\end{eqnarray}

\subsection{Sums of $a_t(p)$}

Using the Explicit Formula, we will find that we need to handle
sums like

\be \sum_{t=N}^{2N} a_t^{r_1}(p_1) \cdots a_t^{r_n}(p_n). \ee

We record these results for later use. Define

\begin{eqnarray}
A_{r,\mathcal{F}}(p) = \sum_{t (p)} a_t^r(p).
\end{eqnarray}

\begin{lem}\label{lemusefulapsum}
Let $p_1, \dots, p_n$ be distinct primes and $r_i \geq 1$. Then
\begin{eqnarray}
\sum_{t (p_1\cdots p_n)} \prod_{i=1}^n a_t^{r_i}(p_i) &=&
\prod_{i=1}^n A_{r_i,\mathcal{F}}(p_i).
\end{eqnarray}
\end{lem}

The proof is a straightforward induction, using the fact that
$a_{t+mp}(p) = a_t(p)$.

Lemma \ref{lemusefulapsum} is our best analogue to the Petersson
formula, which is used in \cite{ILS} to obtain large support for
the density functions.

$\frac{A_{1,\mathcal{F}}(p)}{p}$ is bounded independent of $p$
(\cite{De}). Rosen and Silverman \cite{RSi} proved the following
conjecture of Nagao \cite{Na}:

\begin{thm}[Rosen-Silverman]\label{thmsilvermanrosen} For a one-parameter
family $\mathcal{E}$ of elliptic curves over $\Q(t)$, if Tate's
conjecture is true, then
\end{thm}
\begin{eqnarray}
\lim_{X \to \infty} \frac{1}{X} \sum_{p \leq X}
-\frac{A_{1,\mathcal{F}}(p)}{p} \log p = \mbox{rank} \
\mathcal{E}(\Q(t))
\end{eqnarray}

Tate's conjecture is known for rational surfaces (see \cite{RSi}).
An elliptic surface $y^2 = x^3 + A(t)x + B(t)$ is rational iff one
of the following is true: $(1)\ $ $0 < \max\{3 \mbox{deg} A,
2\mbox{deg} B\} < 12;$ $(2)\ $ $3\mbox{deg} A = 2\mbox{deg} B =
12$ and $\mbox{ord}_{t=0}t^{12} \Delta(t^{-1}) = 0$.

\begin{thm}[Michel \cite{Mi}]\label{thmmichel}
Consider a one-parameter family over $\Q(t)$ with non-constant
$j(E_t)$. Then
\begin{eqnarray}
A_{2,\mathcal{F}}(p) = p^2 + O(p^{\frac{3}{2}}).
\end{eqnarray}
\end{thm}

\subsection{Sieving and Conductors}

To evaluate the sums of $\prod_i a_t^{r_i}(p_i)$, it is necessary
to restrict $t$ to arithmetic progressions; in order to bound some
of the error terms, we will see that the conductors $C(t)$ must be
monotone.

Let

\begin{eqnarray}
\mathcal{T}_{sqfree} &=& \Big\{t \in [N,2N]:\ D(t) \ \mbox{is
sqfree} \Big\} \nonumber\\ \mathcal{T}_N &=& \Big\{t \in [N,2N]:\
d^2 \nmid D(t) \ \mbox{for} \ 2 \leq d \leq \log^l N \Big\}.
\end{eqnarray}

Clearly $\mathcal{T}_{sqfree} \subset \mathcal{T}_N$. We show
$\mathcal{T}_N$ is a union of arithmetic progressions, and
$|\mathcal{T}_N - \mathcal{T}_{sqfree}| = o(N)$.

Thus, except for $o(N)$ values of $t$, we can write $t$ good
(where the conductors are monotone) as a union of arithmetic
progressions. For proofs, see Theorems \ref{thmconditionsnfbig}
and \ref{thmcondcard}.


\section{$1$- and $2$-Level Density Kernels for the Classical Compact
Groups}\label{secksdensities} \setcounter{equation}{0}

By \cite{KS1}, the $n$-level densities for the classical compact
groups are

\begin{eqnarray}\label{eqdensitykernels}
W_{n,O^+}(x) & = & \textbf{det} (K_1(x_i,x_j))_{i,j\leq n}
\nonumber\\ W_{n,O^-}(x) & = & \textbf{det}
(K_{-1}(x_i,x_j))_{i,j\leq n}  + \sum_{k=1}^n \delta(x_k)
\textbf{det}(K_{-1}(x_i,x_j))_{i,j\neq k} \nonumber\\ & = &
(W_{n,O^-})_{1}(x) + (W_{n,O^-})_{2}(x) \nonumber\\ W_{n,O}(x) & =
& \foh W_{n,O^+}(x) + \foh W_{n,O^-}(x) \nonumber\\ W_{n,U}(x) & =
& \textbf{det} (K_0(x_i,x_j))_{i,j\leq n} \nonumber\\ W_{n,Sp}(x)
&=& \textbf{det} (K_{-1}(x_i,x_j))_{i,j\leq n}
\end{eqnarray}

where $K(y) = \kkot{y}$, $ K_\epsilon(x,y) = K(x-y) + \epsilon
K(x+y)$ for $\epsilon = 0, \pm 1$, $O^+$ denotes the group
$SO(\mbox{even})$ and $O^-$ the group $SO(\mbox{odd})$.

\subsection{$1$-Level Densities}

Let $I(u)$ be the characteristic function of $[-1,1]$.

\begin{thm}[$1$-Level Densities]\label{thmonelevel}
\begin{eqnarray}
\widehat{W_{1,O^+} }(u) & = & \delta(u) + \foh I(u) \nonumber\\
\widehat{W_{1,O} }(u) & = & \delta(u) + \foh \nonumber\\
\widehat{W_{1,O^-} }(u) & = & \delta(u) - \foh I(u) + 1
\nonumber\\ \widehat{W_{1,Sp} }(u) & = & \delta(u) - \foh I(u)
\nonumber\\ \widehat{W_{1,U} }(u) & = & \delta(u). \end{eqnarray}
For functions whose Fourier Transforms are supported in $[-1,1]$,
the three orthogonal densities are indistinguishable, though they
are distinguishable from $U$ and $Sp$. To detect differences
between the orthogonal groups using the $1$-level density, one
needs to work with functions whose Fourier Transforms are
supported beyond $[-1,1]$.
\end{thm}

\subsection{$2$-Level Densities}

\begin{thm}[$\mathcal{G} = SO(\mbox{even})$, $O$,
or $SO(\mbox{odd})$]\label{thmtwolevel}

Let $c(\mathcal{G}) = 0$, $\foh$, $1$ for $\mathcal{G} =
SO(\mbox{even})$, $O$, $SO(\mbox{odd})$. For even functions
supported in $|u_1| + |u_2| < 1$
\begin{eqnarray}
& &  \ \ \int \int \widehat{f_1}(u_1)\widehat{f_2}(u_2)
\widehat{W_{2,\mathcal{G}}}(u) du_1du_2 \nonumber\\ & = &
\Big[\hfo(0) + \foh f_1(0) \Big] \Big[\hft(0) + \foh f_2(0) \Big]
\ + \ 2 \int |u| \hfo(u) \hft(u)du \nonumber\\ & & \ - \ 2
\widehat{f_1f_2}(0) \ - \ f_1(0)f_2(0)  \ + \
c(\mathcal{G})f_1(0)f_2(0). \end{eqnarray}

For arbitrarily small support, the three $2$-level densities
differ. One increases by a factor of $\foh f_1(0) f_2(0)$ moving
from $\widehat{W_{2,O^+}}$ to $\widehat{W_{2,O}}$ to
$\widehat{W_{2,O^-}}$.
\end{thm}

\begin{thm}[$\mathcal{G} = Sp$]
\begin{eqnarray}
& & \ \ \int \int \widehat{f_1}(u_1)\widehat{f_2}(u_2)
\widehat{W_{2,Sp}}(u) du_1du_2 \nonumber\\ &= & \Big[\hfo(0) +
\foh f_1(0) \Big] \Big[\hft(0) + \foh f_2(0) \Big] \ + \ 2 \int
|u| \hfo(u) \hft(u)du \nonumber\\ & & - 2 \widehat{f_1f_2}(0) -
f_1(0)f_2(0) - f_1(0)\hft(0) - \hfo(0)f_2(0) + 2f_1(0)f_2(0).
\nonumber\\
\end{eqnarray}
\end{thm}

\begin{thm}[$\mathcal{G} = U$]
\begin{eqnarray}
\int \int \widehat{f_1}(u_1)\widehat{f_2}(u_2) \widehat{W_{2,U}}
du_1du_2 &= & \hfo(0)\hft(0) + \int |u| \hfo(u) \hft(u)du
 - \widehat{f_1f_2}(0). \nonumber\\ \end{eqnarray}
\end{thm}

For test functions with arbitrarily small support, the $2$-level
densities for the classical compact groups are mutually
distinguishable.


\section{Expansions for the $1$- and $2$-Level Densities for Elliptic Curve Families}
\setcounter{equation}{0}

For $i = 1$ and $2$, let $f_i$ be an even Schwartz function whose
Fourier Transform is supported in $(-\sigma_i,\sigma_i)$ and
$f(x_1, x_2) = f_1(x_1)f_2(x_2)$, $\widehat{f}(u_1,u_2) $ $=
\hfo(u_1)\hft(u_2)$.

\subsection{$1$-Level Density: $D_{1,\mathcal{F}}(f)$}

\begin{eqnarray}
D_{1,\mathcal{F}}(f) &=& \oof \sumef \sum_{\gamma_E^{(j)}}
f_1\Big(\gamma_E^{(j)}\frac{\log N_E}{2\pi}\Big) \nonumber\\ & = &
\widehat{f_1}(0) + f_1(0) -  2 \sum_p \oop \oof \sumef \plogne
\widehat{f_1} \Big( \plogne \Big) a_E(p) \nonumber\\ & & - 2
\sum_p \frac{1}{p^2} \oof \sumef \plogne \widehat{f_1}
\Big(\frac{2 \log p}{\log N_E} \Big) a_E^2(p) \nonumber\\ & & +
O\Big(\frac{\log \log N_E}{\log N_E}\Big).
\end{eqnarray}

As the $1$-level density sums are sub-calculations which arise in
the $2$-level investigations, we postpone their determination for
now.

\subsection{$2$-Level Density: $D_{2,\mathcal{F}}(f)$ and
$D_{2,\mathcal{F}}^*(f)$}

Recall the $2$-level density $D_{2,\mathcal{F}}(f)$ is the sum
over all indices $j_1$, $j_2$ with $j_1 \neq \pm j_2$.

\begin{defi}$D_{2,\mathcal{F}}^*(f)$ differs
from the $2$-level density $D_{2,\mathcal{F}}(f)$ in that $j_1$
may equal $\pm j_2$.
\end{defi}

We first calculate $D_{2,\mathcal{F}}^*(f)$, and then subtract off
the contribution from $j_1 = \pm j_2$. Assuming GRH, we may write
the zeros as $1 + i\gamma^{(j)}$, with $\gamma^{(j)} =
-\gamma^{(-j)}$.

\begin{eqnarray}
D_{2,\mathcal{F}}^*(f) & = & \oof \sumef \sum_{j_1} \sum_{j_2}
f_1(L \gamma_E^{(j_1)}) f_2(L \gamma_E^{(j_2)}) \nonumber\\ & = &
\oof \sumef \prod_{i=1}^2 \Bigg[ \hfi(0) + f_i(0) - 2 \sum_{p_i}
\pilogne \oopi \hfi \Big( \pilogne \Big) a_E(p_i) \nonumber\\ & &
\ \ \ \ -2 \sum_{p_i} \pilogne \frac{1}{p_i^2} \hfi \Big(
2\pilogne \Big) a_E^2(p_i) + O\Big(\frac{\log \log N_E}{\log N_E}
\Big) \Bigg] \nonumber\\ & = & \oof \sumef \prod_{i=1}^2 \Bigg[
\hfi(0) + f_i(0) + S_{i,1} + S_{i,2} \Bigg].
\end{eqnarray}

We use Theorem \ref{thmerror} to drop the error terms, as they do
not contribute in the limit as $|\mathcal{F}| \to \infty$. The
astute reader will notice Theorem \ref{thmerror} requires us to
know the $1$-level density, and we have postponed that
calculation; however, in the process of calculating the $2$-level
density we will determine the needed sums for the $1$-level
density (without using Theorem \ref{thmerror} to evaluate them).
Thus, there is no harm in removing the error terms.

There are five types of sums we need to investigate: $\oof \sumef$
$S_{i,1}$, $\oof \sumef$ $S_{i,2}$, $\oof \sumef S_{1,1} S_{2,1}$,
$\oof \sumef S_{1,2}S_{2,2}$, and $\oof \sumef S_{1,1} S_{2,2}$
($i \neq j$). In $S_{i,j}$, $i$ refers to which prime ($p_1$ or
$p_2$), and $j$ the power of $a_E(p_\alpha)$ ($1$ or $2$). The
first and the second are what we need to calculate the one-level
densities.

\subsubsection{$j_1 = \pm j_2$}

Let $\rho = 1 + i\gamma_E^{(j)}$ be a zero. For a curve with even
functional equation, we may label the zeros by

\begin{eqnarray}
\cdots \leq \gamma_E^{(-2)} \leq \gamma_E^{(-1)} \leq 0 \leq
\gamma_E^{(1)} \leq \gamma_E^{(2)} \leq \cdots, \gamma_E^{(-k)} =
- \gamma_E^{(k)}, \end{eqnarray}

while for a curve with odd functional equation we label the zeros
by

\begin{eqnarray}
\cdots \leq \gamma_E^{(-1)} \leq 0 \leq \gamma_E^{(0)} = 0 \leq
\gamma_E^{(1)} \leq \cdots, \gamma_E^{(-k)} = - \gamma_E^{(k)}.
\end{eqnarray}

We exclude the contribution from $j_1 = \pm j_2$. If an elliptic
curve has even functional equation, $j_i$ ranges over all non-zero
integers, and $\gamma_E^{(-j)} = -\gamma_E^{(j)}, j \neq -j$.
Since the test functions are even, the sum over all pairs
$(j_1,j_2)$ with $j_1 = \pm j_2$ is twice the sum over all pairs
$(j,j)$, which is $D_{1,E}(f_1f_2)$, ie, the $1$-level density for
a curve $E$ with test function $f_1(x)f_2(x)$.

If an elliptic curve has odd functional equation, $j_i$ ranges
over all integers. The curve vanishes to odd order at the critical
point $s = 1$. Except for one zero (labelled $\gamma_E^{(0)}$),
for every non-zero $j$, $\gamma_E^{(-j)} = - \gamma_E^{(j)}$, and
$j \neq -j$. Twice the sum over pairs $(j,j)$ minus the
contribution from the pair $(0,0)$ equals the sum over all pairs
$(j_1,j_2)$ with $j_1 = \pm j_2$. Thus, the curves with odd sign
contribute $D_{1,E}(f_1f_2) - f_1(0)f_2(0)$.

Let $\epsilon_E = \pm 1$ be the sign of the functional equation
for $E$, and define

\begin{defi} $N(\mathcal{F},-1) = \oof \sumef \frac{1-\epsilon_E}{2}$,
ie, the percent of curves with odd sign.
\end{defi}

Summing over $E \in \mathcal{F}$ yields $D_{1,\mathcal{F}}(f_1f_2)
- N(\mathcal{F},-1)f_1(0)f_2(0)$ for $j_1 = \pm j_2$.

\subsubsection{$2$-Level Density Expansion}

\begin{lem}[$2$-Level Density Expansion]\label{lemdtwoexpansion}
\begin{eqnarray}
D_{2,\mathcal{F}}(f) &=& \oof \sumef \prod_{i=1}^2 \Bigg[ \hfi(0)
+ f_i(0) + S_{i,1} + S_{i,2} \Bigg] \nonumber\\ & & - \
2D_{1,\mathcal{F}}(f_1f_2) + (f_1f_2)(0)N(\mathcal{F},-1) + O\Big(
\frac{\log \log N}{\log N}\Big).\nonumber\\ \end{eqnarray}

To evaluate the above, we only need to know the percent of curves
with odd sign, not which curves are even or odd. For the $3$ and
higher level densities, we have to execute sums over the subset of
curves with odd sign.
\end{lem}

\subsection{Useful Expansion for the $1$- and $2$-Level Densities for
One Parameter Families}

Let $\mathcal{E}$ denote a one-parameter family of elliptic curves
$E_t$ over $\Qoft$, $t \in [N,2N]$, and $\mathcal{F}$ denote a
sub-family of $\mathcal{E}$. In the applications, $\mathcal{F}$
will be obtained by sieving to $D(t)$ good, where $D(t)$ is the
product of the irreducible polynomial factors of $\Delta(t)$.

\subsubsection{Needed Prime Sums}

\begin{lem}[Prime Sums]\label{lemprimesums} Let $C(N)$ be a power of $N$.
By Lemmas \ref{primeonesuma}, \ref{primeonesumb} and
\ref{primetwosuma},
\begin{enumerate}
\item $\sum_p \frac{\log p}{\log C(N)} \oop \hfo\Big( \frac{\log p}{\log C(N)}
\Big) = \foh f_1(0) + O\Big(\frac{1}{\log N}\Big)$
\item $\sum_p \frac{\log p}{\log C(N)} \oop \hfo\Big(2 \frac{\log p}{\log C(N)}
\Big) = \frac{1}{4} f_1(0) + O\Big(\frac{1}{\log N}\Big)$
\item $\sum_p \frac{\log^2 p}{\log^2 C(N)} \frac{1}{p} \hfo
\hft\Big(\frac{\log p}{\log C(N)} \Big) = \foh
\int_{-\infty}^{\infty} |u| \hfo(u) \hft(u) du +
O\Big(\frac{1}{\log N}\Big)$
\end{enumerate}
If instead we are summing over primes congruent to $a$ mod $m$, we
use Lemma \ref{lemprimeonesum} and \ref{primetwosumb}, and the
right-hand sides are modified by $\frac{1}{\varphi(m)}$.
\end{lem}

\subsubsection{Expansions of Sums}

We use the expansion from Lemma \ref{lemdtwoexpansion}. Recall

\begin{eqnarray}
S_{i,j} & = & -2\sum_{p_i} \frac{\log p_i}{\log C(t)}
\frac{1}{p_i^j} \hfi\Big(2^{j-1} \frac{\log p_i}{\log C(t)} \Big)
a_t^{j}(p_i).
\end{eqnarray}

In $S_{i,j}$, $i$ refers to the prime ($p_1, p_2$) and $j$ refers
to the power of $a_t(p)$ ($a_t(p), a_t^2(p)$).

To determine the $1$- and $2$-level densities, there are eight
sums over $t \in \mathcal{F}$ to evaluate: $\oof \sumtinF S_{1,1}$
and $\oof \sumtinF S_{2,1}$; $\oof \sumtinF S_{1,2}$ and $\oof
\sumtinF S_{2,2}$; $\oof \sumtinF S_{1,1} S_{2,2}$ and $\oof
\sumtinF S_{2,1} S_{1,2}$; $\oof \sumtinF S_{1,1} S_{2,1}$; $\oof
\sumtinF S_{1,2} S_{2,2}$.

We have written the sums in pairs where the two sums are handled
similarly. Substituting the definitions leads to five types of
sums:

\begin{enumerate}
\item $-2\sum_p \oop \ooNF \sumtinF \frac{\log p}{\log
C(t)} \hfo\Big(\frac{\log p}{\log C(t)} \Big) a_t(p)$
\item $-2\sum_p \frac{1}{p^2} \ooNF \sumtinF \frac{\log p}{\log
C(t)} \hfo\Big(2\frac{\log p}{\log C(t)} \Big) a_t^2(p)$
\item $4\sum_{p_1} \sum_{p_2} \frac{1}{p_1p_2^2} \ooNF \sumtinF
\frac{\log p_1}{\log C(t)} \frac{\log p_2}{\log C(t)}
\hfo\Big(\frac{\log p}{\log C(t)} \Big)\hft\Big(2\frac{\log
p}{\log C(t)} \Big) a_t(p_1)a_t^2(p_2)$
\item $4\sum_{p_1} \sum_{p_2} \frac{1}{p_1p_2} \ooNF \sumtinF
\frac{\log p_1}{\log C(t)} \frac{\log p_2}{\log C(t)}
\hfo\Big(\frac{\log p}{\log C(t)} \Big)\hft\Big(\frac{\log p}{\log
C(t)} \Big) a_t(p_1)a_t(p_2)$
\item $4\sum_{p_1} \sum_{p_2} \frac{1}{p_1^2p_2^2} \ooNF \sumtinF
\frac{\log p_1}{\log C(t)} \frac{\log p_2}{\log C(t)}
\hfo\Big(2\frac{\log p}{\log C(t)} \Big)\hft\Big(2\frac{\log
p}{\log C(t)} \Big) a_t^2(p_1)a_t^2(p_2)$
\end{enumerate}

In the above sums, we use Lemma \ref{lemsmallprimes} to restrict
to primes greater than $\log^l N$, $l < 2$. Label the five sums
$\oof \sum_{t \in \mathcal{F}} S(t;p)$  by $T_k(p)$ and
$T_k(p_1,p_2)$. Trivially by Hasse some of the above do not
contribute.

In the third sum, if $p_1 = p_2 = p$, we get $\ll \frac{1}{\log N}
\sum_p$ $\frac{p^{\frac{3}{2}}\log p}{p^3}$ $= O(\frac{1}{\log
N})$. In the fifth sum, if $p_1 = p_2 = p$ we get $\ll
\frac{1}{\log N} \sum_p$ $\frac{p^2\log p}{p^4}$ $=
O(\frac{1}{\log N})$.

Thus, we only study the third and fifth sums when $p_1 \neq p_2$.
The fourth sum has the potential to contribute when $p_1 = p_2$.
Hence we break it into two cases: $p_1 \neq p_2$ and $p_1 = p_2$.

\subsubsection{Conditions on the Family to Evaluate the Sums}

\begin{eqnarray}\label{eqconditionsonF}\label{eqfivesums}
 \mbox{Conditions on the
Family $\mathcal{F}$}
\end{eqnarray}
Let $T_k(p)$ and $T_k(p_1,p_2)$ ($= \oof \sum_{t \in \mathcal{F}}
S(t;p)$ ) equal
\begin{enumerate}
\item $\frac{\log p}{\log C(N)} \hfo\Big(\frac{\log p}{\log C(N)} \Big)
\Bigg[-r + O\Big(p^{-\alpha} + \frac{p^\beta}{|\mathcal{F}|} +
\frac{1}{\log^\gamma N} \Big)\Bigg]$
\item $\frac{\log p}{\log C(N)} \hfo\Big(2\frac{\log p}{\log C(N)} \Big)
\Bigg[p + O\Big(p^{1-\alpha} + \frac{p^\beta}{|\mathcal{F}|} +
\frac{p}{\log^\gamma N} \Big)\Bigg]$
\item $\frac{\log p_1}{\log C(N)} \frac{\log p_2}{\log
C(N)} \hfo\Big(\frac{\log p_1}{\log C(N)} \Big)
\hft\Big(2\frac{\log p_2}{\log C(N)} \Big) \Bigg[ -rp_2 +
O\Big(p_1^{-\alpha_1}p_2^{1-\alpha_2} + \frac{
p_1^{\beta_1}p_2^{\beta_2}  }{|\mathcal{F}|} +
\frac{p_2}{\log^\gamma N} \Big) \Bigg]$
\item
\begin{enumerate}
\item $ \frac{\log p_1}{\log C(N)} \frac{\log p_2}{\log
C(N)} \hfo\Big(\frac{\log p_1}{\log C(N)} \Big)
\hft\Big(\frac{\log p_2}{\log C(N)} \Big) \Bigg[ r^2 +
O\Big(p_1^{1-\alpha_1}p_2^{1-\alpha_2} + \frac{
p_1^{\beta_1}p_2^{\beta_2} }{|\mathcal{F}|} + \frac{1}{\log^\gamma
N} \Big) \Bigg]$  if $p_1 \neq p_2$
\item $\frac{\log^2 p}{\log^2 C(N)} \hfo\hft\Big(\frac{\log p}{\log C(N)}
\Big) \Bigg[ p  + O\Big(p^{1-\alpha} +
\frac{p^\beta}{|\mathcal{F}|} + \frac{p}{\log^\gamma N} \Big)
\Bigg]$ if $p_1 = p_2 = p$
\end{enumerate}
\item $\frac{\log p_1}{\log C(N)} \frac{\log p_2}{\log
C(N)} \hfo\Big(2\frac{\log p_1}{\log C(N)} \Big)
\hfo\Big(2\frac{\log p_2}{\log C(N)} \Big) \Bigg[ p_1p_2  +
O\Big(p_1^{1-\alpha_1}p_2^{1-\alpha_2} + \frac{
p_1^{\beta_1}p_2^{\beta_2} }{|\mathcal{F}|} + \frac{p_1
p_2}{\log^\gamma N} \Big) \Bigg]$
\end{enumerate}
where $\alpha, \beta, \gamma > 0$, $\alpha_i, \beta_i \geq 0$ and
whenever two $\alpha_i$ or $\beta_i$ occur, at least one is
positive.

By Lemma \ref{lemprimesums} we can evaluate the eight $S_{i,j}$
sums for a family satisfying Conditions \ref{eqconditionsonF}:

\begin{lem}[$S_{i,j}$ Sums]\label{lemvaluesofSij} If the family
satisfies Conditions \ref{eqconditionsonF}, then (up to lower
order terms which do not contribute for small support),
\begin{enumerate}
\item $\oof \sumtinF S_{i,1} = rf_i(0)$
\item $\oof \sumtinF S_{i,2} = -\foh f_i(0) $
\item $\oof \sumtinF S_{1,1} S_{2,2} + S_{2,1} S_{1,2} = -\foh r f_1(0)f_2(0) +
-\foh rf_1(0)f_2(0)$
\item $\oof \sumtinF S_{1,1} S_{2,1} = r^2 f_1(0)f_2(0) +
2\int_{-\infty}^\infty |u|\hfo(u)\hft(u)du$
\item $\oof \sumtinF S_{1,2} S_{2,2} = \frac{1}{4}f_1(0)f_2(0)$
\end{enumerate}
\end{lem}

\subsubsection{$1$- and $2$-Level Densities, Assuming Certain
Conditions on the Family}

Substituting Lemma \ref{lemvaluesofSij} into the $1$- and
$2$-level density expansions we obtain

\begin{lem}[$1$- and $2$-Level
Densities]\label{lemrankwithassumptions} Assume $\NF$ is a
positive multiple of $N$ and $\mathcal{F}$ satisfies conditions
\ref{eqconditionsonF}. Up to lower order correction terms (which
vanish as $\NF \to \infty$), for even Schwartz functions with
small support,
\begin{eqnarray}
D_{1,\mathcal{F}}(f) = \hfo(0) + \foh f_1(0) + rf_1(0)
\end{eqnarray} and
\begin{eqnarray}
D_{2,\mathcal{F}}(f) &=& \prod_{i=1}^2 \Bigg[ \hfi(0) + \foh
f_i(0) \Bigg] + 2\int_{-\infty}^\infty |u|\hfo(u)\hft(u)du
\nonumber\\ & &  - 2\widehat{f_1f_2}(0) - f_1(0)f_2(0) +
(f_1f_2)(0)N(\mathcal{F},-1) \nonumber\\ & & + (r^2-r)f_1(0)f_2(0)
+ r\hfo(0)f_2(0) + r f_1(0)\hft(0). \end{eqnarray}

\textbf{Let $D_{1,\mathcal{F}}^{(r)}(f_1)$ and
$D_{2,\mathcal{F}}^{(r)}(f_1)$ be the $1$- and $2$-level densities
from which the contributions of $r$ family zeros at the critical
point have been subtracted}. Then
\begin{eqnarray}
D_{1,\mathcal{F}}^{(r)}(f_1) = \hfo(0) + \foh f_1(0)
\end{eqnarray} and
\begin{eqnarray}
D_{2,\mathcal{F}}^{(r)}(f_1) &=& \prod_{i=1}^2 \Bigg[ \hfi(0) +
\foh f_i(0) \Bigg] + 2\int_{-\infty}^\infty |u|\hfo(u)\hft(u)du
\nonumber\\ & &  - 2\widehat{f_1f_2}(0) - f_1(0)f_2(0) +
(f_1f_2)(0)N(\mathcal{F},-1). \end{eqnarray} Thus, removing the
contribution from $r$ family zeros, for test functions of small
support the $2$-level density of the remaining zeros agrees with
$SO(\mbox{even})$ if all curves are even, $O$ if half are even and
half odd, and $SO(\mbox{odd})$ if all are odd.
\end{lem}

Proof: The $1$-level density is immediate from substitution.
Substituting for the eight $S_{i,j}$ sums for
$D_{2,\mathcal{F}}(f)$ yields (up to lower order terms which don't
contribute for small support)

\begin{eqnarray}
D_{2,\mathcal{F}}(f) &=& = \prod_{i=1}^2 \Bigg[ \hfi(0) + f_i(0)
\Bigg] \nonumber\\ & & + \ \Bigg[ \hfo(0) + f_1(0) \Bigg] rf_2(0)
+ \Bigg[ \hft(0) + f_2(0) \Bigg] rf_1(0) \nonumber\\ & &  + r^2
f_1(0)f_2(0) + 2\int_{-\infty}^\infty |u|\hfo(u)\hft(u)du
\nonumber\\ & & + \ \Bigg[ \hfo(0) + f_1(0) \Bigg] \Big(-\foh
f_2(0)\Big) + \Bigg[ \hft(0) + f_2(0) \Bigg] \Big(-\foh
f_1(0)\Big) \nonumber\\ & & - \foh r f_1(0)f_2(0) - \foh
rf_1(0)f_2(0) + \frac{1}{4}f_1(0)f_2(0) \nonumber\\ & & - \
2D_{1,\mathcal{F}}(f_1f_2) + (f_1f_2)(0)N(\mathcal{F},-1) + O\Big(
\frac{\log \log N}{\log N}\Big) \nonumber\\ & =& \prod_{i=1}^2
\Bigg[ \hfi(0) + \foh f_i(0) \Bigg] + 2\int_{-\infty}^\infty
|u|\hfo(u)\hft(u)du \nonumber\\ & & + 2rf_1(0)f_2(0) +
r\hfo(0)f_2(0) + r f_1(0)\hft(0) - rf_1(0)f_2(0) + r^2
f_1(0)f_2(0) \nonumber\\ & & -
 2D_{1,\mathcal{F}}(f_1f_2) + (f_1f_2)(0)N(\mathcal{F},-1). \end{eqnarray}

Substituting

\begin{eqnarray}
D_{1,\mathcal{F}}(f_1f_2) = \widehat{f_1f_2}(0) + \foh
f_1(0)f_2(0) + rf_1(0)f_2(0) \end{eqnarray}

yields

\begin{eqnarray}
D_{2,\mathcal{F}}(f) &=& \prod_{i=1}^2 \Bigg[ \hfi(0) + \foh
f_i(0) \Bigg] + 2\int_{-\infty}^\infty |u|\hfo(u)\hft(u)du
\nonumber\\ & & + rf_1(0)f_2(0) + r\hfo(0)f_2(0) + r f_1(0)\hft(0)
+ r^2 f_1(0)f_2(0) \nonumber\\ & & - 2\widehat{f_1f_2}(0) -
f_1(0)f_2(0) - 2rf_1(0)f_2(0) + (f_1f_2)(0)N(\mathcal{F},-1)
\nonumber\\ & = & \prod_{i=1}^2 \Bigg[ \hfi(0) + \foh f_i(0)
\Bigg] + 2\int_{-\infty}^\infty |u|\hfo(u)\hft(u)du \nonumber\\ &
&  - 2\widehat{f_1f_2}(0) - f_1(0)f_2(0) +
(f_1f_2)(0)N(\mathcal{F},-1) \nonumber\\ & & + (r^2-r)f_1(0)f_2(0)
+ r\hfo(0)f_2(0) + r f_1(0)\hft(0). \end{eqnarray}

If the family has rank $r$ over $\Qoft$, there is a natural
interpretation of these terms. By the Birch and Swinnerton-Dyer
conjecture (used only for interpretation purposes) and Silverman's
Specialization Theorem, for all $t$ sufficiently large, each
curve's $L$-function has at least $r$ zeros at the critical point.
We isolate the contributions from $r$ family zeros.

Assume there are $r$ family zeros at the critical point. Let $L_t
= \frac{\log C(t)}{2\pi}$. Recall the $1$-level density is
$D_{1,\mathcal{F}}(f) = \widehat{f}(0) + \foh f(0) + rf(0)$. Let
$j_i$ range over all zeros of a curve, and $j_i'$ range over all
but the $r$ family zeros.

\begin{eqnarray}
D_{2,\mathcal{F}}(f) &=& \oof \sumtinF \sum_{j_1} \sum_{j_2}
f_1(L_t\gamma_{E_t}^{(j_1)}) f_2(L_t\gamma_{E_t}^{(j_2)})
\nonumber\\ & & - \ 2D_{1,\mathcal{F}}(f_1f_2) +
(f_1f_2)(0)N(\mathcal{F},-1) \nonumber\\ &=& \oof \sumtinF
\Big(rf_1(0) + \sum_{j_1'} f_1(L_t\gamma_{E_t}^{(j_1')}) \Big)
\Big(rf_2(0) + \sum_{j_2'} f_2(L_t\gamma_{E_t}^{(j_2')}) \Big)
\nonumber\\ & & - \ 2D_{1,\mathcal{F}}(f_1f_2) +
(f_1f_2)(0)N(\mathcal{F},-1) \nonumber\\ &=& \oof \sumtinF
\sum_{j_1'} \sum_{j_2'} f_1(L_t\gamma_{E_t}^{(j_1')})
f_2(L_t\gamma_{E_t}^{(j_2')}) \nonumber\\ & & + rf_1(0)
D_{1,\mathcal{F}}(f_2) + D_{1,\mathcal{F}}(f_1) rf_2(0) -
r^2f_1(0)f_2(0) \nonumber\\ & & - 2D_{1,\mathcal{F}}(f_1f_2) +
(f_1f_2)(0)N(\mathcal{F},-1) \nonumber\\ &=& \oof \sumtinF
\sum_{j_1'} \sum_{j_2'} f_1(L_t\gamma_{E_t}^{(j_1')})
f_2(L_t\gamma_{E_t}^{(j_2')}) + (f_1f_2)(0)N(\mathcal{F},-1)
\nonumber\\ & & + rf_1(0)\Big(\hft(0) + (r+\foh)f_2(0)\Big) +
\Big(\hfo(0) + (r+\foh)f_1(0)\Big)rf_2(0) \nonumber\\ & & -
r^2f_1(0)f_2(0) - 2\Big(\widehat{f_1f_2}(0) + \foh f_1(0)f_2(0) +
rf_1(0)f_2(0)\Big) \nonumber\\ &=& \Bigg[ \oof \sumtinF
\sum_{j_1'} \sum_{j_2'} f_1(L_t\gamma_{E_t}^{(j_1')})
f_2(L_t\gamma_{E_t}^{(j_2')}) \nonumber\\ & & - 2\Big(
\widehat{f_1f_2}(0) + \foh f_1(0)f_2(0) \Big) +
(f_1f_2)(0)N(\mathcal{F},-1)\Bigg] \nonumber\\ & & +
rf_1(0)\hft(0) + r\hfo(0)f_2(0) + (r^2-r)f_1(0)f_2(0) \nonumber\\
&=& D_{2,\mathcal{F}}^{(r)}(f_1) + rf_1(0)\hft(0) + r\hfo(0)f_2(0)
+ (r^2-r)f_1(0)f_2(0). \end{eqnarray}

We isolate

\begin{lem} The contribution from $r$ critical point zeros is
\begin{eqnarray}
rf_1(0)\hft(0) + r\hfo(0)f_2(0) + (r^2-r)f_1(0)f_2(0).
\end{eqnarray}
\end{lem}


\section{Calculation of the $1$- and $2$-Level
Densities for Elliptic Curve Families} \setcounter{equation}{0}

Let $\mathcal{E}$ be a one-parameter family of elliptic curves
$E_t$ with discriminants $\Delta(t)$ and conductors $C(t)$. For
many families, we can evaluate the conductors exactly if we sieve
to a subfamily $\mathcal{F}$ defined as the $t \in [N,2N]$ with
$D(t)$ good, where $D(t) = a_k t^k + \cdots + a_0$ ($a_k \geq 1$)
is the product of the irreducible polynomial factors of
$\Delta(t)$. Usually good will mean square-free, although
occasionally it will mean square-free except for a fixed set of
primes, and for these special primes, the power of $p|D(t)$ is
independent of $t$.

Let our family $\mathcal{F}$ be the set of good $t \in [N,2N]$
where the conductors are given by a monotone polynomial in $t$. We
use this polynomial for the conductors at non-good $t$; this is
permissible as these curves are not in our family, and do not
originally appear in our sums.

For each $d$, let

\begin{eqnarray}
T(d) = \{t \in [N,2N]: d^2|D(t)\}. \end{eqnarray}

Let $S(t)$ be some quantity associated to the elliptic curve
$E_t$. We study

\begin{eqnarray}
\sum_{t=N \atop D(t) \ good }^{2N} S(t) =
\sum_{d=1}^{(2a_kN)^{\frac{k}{2}}} \mu(d) \sum_{t \in T(d)} S(t).
\end{eqnarray}

In particular, setting $S(t) = 1$ yields the cardinality of the
family. In all the families we investigate, $|\mathcal{F}| =
c_{\mathcal{F}}N + o(N)$, $c_{\mathcal{F}} > 0$.

Let $t_1(d), \dots, t_{\nu(d)}(d)$ be the incongruent roots of
$D(t) \equiv 0$ mod $d^2$. The presence of $\mu(d)$ allows us to
restrict to $d$ square-free. For small $d$, we may take the
$t_i(d) \in [N,N + d^2)$. For such $d$,

\begin{eqnarray}
\sum_{t \in T(d)} S(t) = \sum_{i=1}^{\nu(d)} \sum_{t' =
0}^{[N/d^2]} S\Big(t_i(d) + t'd^2\Big) \ + \
O\Big(\nu(d)||S||_\infty \Big).
\end{eqnarray}

The error piece is from boundary effects for the last value of
$t'$. $T(d)$ restricts us to $t \in [N,2N]$; as each $t_i(d) \geq
N$, and at most one is exactly $N$, it is possible in summing to
$t' = [N/d^2]$ we've added an extra term.

\subsection{Assumptions for Sieving}

We evaluate the sums under the following assumptions:

\begin{enumerate}
  \item For square-free $D(t)$, the conductors $C(t)$ are given by
  a monotone polynomial in $t$.
  \item A positive percent of $t \in [N,2N]$ have $D(t)$
  square-free; ie, $|\mathcal{F}| = c_{\mathcal{F}}N + o(N)$.
\end{enumerate}

We constantly use Lemma \ref{lemnudbound} ($\nu(d) \ll d^\epsilon$
for square-free $d$) and

\begin{eqnarray}\label{eqsizeoffamily}
\sum_{t = N \atop D(t) \ good}^{2N} 1 &=& \sum_{d=1}^{\log^l N}
\mu(d) \sum_{t = N \atop D(t) \equiv 0 (d^2)}^{2N} 1 \ + \ o(N) =
c_{\mathcal{F}} N + o(N), \ c_{\mathcal{F}} > 0. \ \ \
\end{eqnarray}

We show the family satisfies Conditions \ref{eqconditionsonF}. We
evaluate the sums over $t \in \mathcal{F}$ below and then execute
the summation over the prime(s). $\hfi$ is supported in
$(-\sigma_i,\sigma_i)$. There are no contributions (for $\sigma_i$
sufficiently small) in the prime sum(s) for sufficiently small
error terms.

\subsection{Definition of Terms for Sieving}

Recall $A_{r,\mathcal{F}}(p) = \sum_{t (p)} a_t^r(p)$. For
distinct primes, by Lemma \ref{lemusefulapsum}

\begin{eqnarray}
\sum_{t (p_1\cdots p_n)} \prod_{j=1}^n a_t^{r_i}(p_j) &=&
\prod_{j=1}^n A_{r_i,\mathcal{F}}(p_i).
\end{eqnarray}

By Lemma \ref{lemsmallprimes}, we may assume all of our primes (in
the expansion from the Explicit Formula in the $n$-level
densities) are at least $\log^l N$, $l \in [1,2)$. We can
incorporate these errors into our existing error terms; the result
will still be a lower order term which will not contribute for
small support.

$S(t)$ will equal $\widetilde{a}_P(t) G_P(t)$, where for distinct
primes $p_1$ and $p_2$

\begin{eqnarray}
\widetilde{a}_P(t) & = & a_t^{r_1}(p_1) a_t^{r_2}(p_2) \nonumber\\
G_P(t) & = & \prod_{j=1 \atop r_j \neq 0}^2 \frac{\log p_j}{\log
C(t)} f_j\Big(2^{r_j-1}\frac{\log p_j}{\log C(t)} \Big)
\nonumber\\ (r_1,r_2) & \in & \Big\{(1,0), (0,1), (2,0), (1,1),
(0,2), (1,2), (2,1), (2,2) \Big\}. \end{eqnarray}

Thus $\widetilde{a}_P(t) G_P(t)$ is merely a convenient way of
encoding the eight sums we need to examine for the $1$ and
$2$-level densities.

Actually, this is slightly off. We have to study

\begin{eqnarray}
\prod_{j=1 \atop r_j \neq 0}^2 \frac{1}{p_j^{r_j}} \frac{\log
p_j}{\log C(t)} g_j\Big(2^{r_j-1} \frac{\log p_j}{\log C(t)} \Big)
a_t^{r_j}(p_j). \end{eqnarray}

If both $r_j$'s are non-zero and the two primes are equal, we
obtain

\begin{eqnarray}
\frac{1}{p^{r_1+r_2}} \Big(\frac{\log p}{\log C(t)}\Big)^2 \times
\cdots \times a_t^{r_1+r_2}(p). \end{eqnarray}

For example, if $r_1 = r_2 = 1$ we would get $(\frac{\log p}{\log
C(t)})^2 \times \cdots \times a_t^2(p)$. Thus, the definition of
$G_P$ needs to be slightly modified. We want to deal with distinct
primes $p_1$ and $p_2$. There will be no contribution for equal
primes if $r_1 + r_2 \geq 3$; simply bound each $a_t(p)$ by Hasse.
There is a contribution if $r_1 = r_2 = 1$. By modifying the
definition of $G_P$ we may regard it as a case where $r = (2,0)$;
however, we have $(\frac{\log p}{\log C(t)})^2$ instead of
$(\frac{\log p}{\log C(t)})$, and instead of $f_1(\cdots)$ we will
have $f_1f_2(\cdots)$. Note we evaluate the test functions at
$\frac{\log p}{\log C(t)}$ and not $2 \frac{\log p}{\log C(t)}$.
We have

\begin{eqnarray}
G_P(t) &=& \prod_{j=1 \atop r_j \neq 0}^2 \Big(\frac{\log
p_j}{\log C(t)}\Big)^{\kappa(r)}
g_j\Big(2^{r_j-\kappa(r)}\frac{\log p_j}{\log C(t)} \Big),
\end{eqnarray}

where $\kappa(r)$ is $2$ if $r = (2,0)$ and this arises from $p_1
= p_2 = p$ and $\kappa(r) = 1$ otherwise; $g_j = f_j$ unless $r =
(2,0)$ arising from $p_1 = p_2 = p$, in which case $g_1 = f_1f_2$.

We may now assume the primes are distinct. Define

\begin{eqnarray}
P &=& \prod_{j=1 \atop r_j \neq 0}^2 p_j, \ \ \ r = (r_1,r_2), \
r_j \in \{0,1,2\} \nonumber\\ S_c(r,P) & = & \sum_{t (P)}
\widetilde{a}_P(t) = \sum_{t (P)} a_t^{r_1}(p_1) a_t^{r_2}(p_2) =
A_{r_1,\mathcal{F}}(p_1) A_{r_2,\mathcal{F}}(p_2), \end{eqnarray}

where for convenience we set $A_0(p) = 1$. We often have
incomplete sums of $\widetilde{a}_P(t)$ mod $P$. Let $S_I(r,P)$
denote a generic incomplete sum. By Hasse,

\begin{eqnarray}
S_I(r,P) & \leq & P \cdot 2^{r_1} \sqrt{p_1^{r_1}} \cdot 2^{r_2}
\sqrt{p_2^{r_2}} = 2^{r_1 + r_2} p_1^{1 + \frac{r_1}{2}} \cdot
p_2^{1 + \frac{r_2}{2}} = 2^r P^{1 + \frac{r}{2}}, \nonumber\\
\end{eqnarray}

where the last expression is a convenient abuse of notation:

\begin{eqnarray}
2^r & = & 2^{r_1 + r_2}, \ \ \ P^r = p_1^{r_1}\cdot p_2^{r_2}.
\end{eqnarray}

For fixed $i$ and $d$, we evaluate the arguments at $t = t_i(d) +
t'd^2$. Let

\begin{eqnarray}
\widetilde{a}_{d,i,P}(t') &=& \widetilde{a}_P\Big(t_i(d) +
t'd^2\Big), \ \ \  G_{d,i,P}(t') = G_P\Big(t_i(d) +
t'd^2\Big).\nonumber\\
\end{eqnarray}

\subsection{Ranges and Contributions of Sums over Primes}

Each prime sum is to (approximately) $C(N)^{\frac{\sigma_j}{2^{r_j
- \kappa(r)}}} \approx N^{\frac{m\sigma_j}{2^{r_j - \kappa(r)}}}$,
as $C(t)$ is a degree $m$ polynomial. We assume $\sigma_j < \foh$
as we do not worry about $p^2 > N$. This is harmless, as handling
the error terms forces the support to be significantly less than
$\foh$.

\begin{lem}[Contributions from Sums over
Primes]\label{lemprimescontributerange} For $r_j = 1$, summing
$\frac{p^\foh}{|\mathcal{F}|}$ does not contribute for $\sigma_j <
\frac{2}{3 m}$. For $r_j = 2$, summing $\frac{1}{|\mathcal{F}|}$
does not contribute for $\sigma_j < \frac{2}{m}$ for $\kappa(r) =
1$ and $\frac{1}{m}$ for $k(r) = 2$. As we often have two sums,
dividing the above supports by $2$ ensures all errors are
manageable: write $\frac{1}{|\mathcal{F}|}$ as
$\frac{1}{\sqrt{|\mathcal{F}|}} \frac{1}{\sqrt{|\mathcal{F}|}}$.
\end{lem}

\subsubsection{Expected Result}

To simplify the proof, we assume

\begin{eqnarray}\label{eqexpectedAisums}
A_{1,\mathcal{F}}(p) & = & -rp + O(1) \nonumber\\
A_{2,\mathcal{F}}(p) & = & p^2 + O(p^{\frac{3}{2}}).
\end{eqnarray}

For a general rational surface, $A_{1,\mathcal{F}}(p) \neq -rp +
O(1)$. A careful book-keeping of the arguments below show that we
only need to be able to handle sums such as

\begin{eqnarray}
\sum_p \frac{\log p}{\log X} f\Big(\frac{\log p}{\log X}\Big)
\frac{A_{1,\mathcal{F}}(p)}{p^2}. \end{eqnarray}

For surfaces where Tate's conjecture is known, we may replace
$A_{1,\mathcal{F}}(p)$ in the above sum with the rank of the
family over $\Qoft$ (see Lemma \ref{lemaepone} and \cite{RSi}).
For notational simplicity, in the proof below we assume
$A_{1,\mathcal{F}}(p) = -rp + O(1)$, and content ourselves with
noting a similar proof works in general.

$A_{r_j}(p_j) = c_j \cdot p_j^{r_j}$ plus lower order terms not
contributing for any support. (This is not quite true. For
families where the curves have complex multiplication, often
$a_t(p)$ vanishes for half the primes, and has double the expected
contribution for the other primes. This case is handled similarly,
using Lemmas \ref{lemprimeonesum} and \ref{primetwosumb}).

Hence $S_c(r,P) = c_1 c_2 p_1^{r_1}p_2^{r_2} = c_1 c_2 P^r$ plus
lower terms. For each pair $(d,i)$ we expect (if we can manage the
conductors) to have approximately $\frac{N/d^2}{P}$ complete sums
of $S_c(r,P) = c_1 c_2 P^r$. We hit this with $\frac{1}{N}
\frac{\log p_j}{\log C(t)} \frac{1}{p_j^{r_j}}$ for each non-zero
$r_j$. We have approximately $\frac{\log p_j}{\log C(t)}
\frac{1}{P^r}$.

A sum like $\sum_{p_j} \frac{\log p_j}{\log C(t)} \frac{1}{p_j}
g(\frac{\log p_j}{\log C(t)})$ contributes; if we had an
additional $\frac{1}{\log N}$ there would be no net contribution.

Thus, we expect terms of the size $P^r$ to contribute, and
$\frac{P^r}{\log N}$ to not contribute.

We rewrite Conditions \ref{eqconditionsonF} in a more tractable
form, using $A_{1,\mathcal{F}}(p)$, $A_{2,\mathcal{F}}(p)$ and
$S_c(r,P)$. Assume the family satisfies Equation
\ref{eqexpectedAisums} (or the related equation if $a_t(p)$
vanishes for half the primes). Then

\begin{enumerate}
\item $P = p$, $\widetilde{a}_P(t) = a_t(p)$: $\frac{S_c(r,P)}{P} =
\frac{-rp + O(1)}{p} = -r + O(\frac{1}{p})$
\item $P = p$, $\widetilde{a}_P(t) = a_t^2(p)$: $\frac{S_c(r,P)}{P} =
\frac{p^2 + O(p^{\frac{3}{2}})}{p} = p + O(\sqrt{p})$
\item $P = p_1p_2$, $\widetilde{a}_P(t) = a_t(p_1)a_t^2(p_2)$:
$\frac{S_c(r,P)}{P} = \frac{-rp_1p_2^2 +
O(p_1p_2^{\frac{3}{2}})}{p_1p_2} = -rp_2 + O(\sqrt{p_2})$
\item $P = p_1p_2$, $\widetilde{a}_P(t) = a_t(p_1)a_t(p_2)$:
\begin{enumerate}
\item $\frac{S_c(r,P)}{P} = \frac{r^2p_1p_2 + O(p_1+p_2)}{p_1p_2} =
r^2 + O(\sqrt{p_1} +\sqrt{p_2})$ if $p_1 \neq p_2$
\item $\frac{S_c(r,P)}{P} = \frac{p^2 + O(p^{\frac{3}{2}})}{p} = p + O(\sqrt{p})$
if $p_1 = p_2 = p$
\end{enumerate}
\item $P = p_1p_2$, $\widetilde{a}_P(t) = a_t^2(p_1)a_t^2(p_2)$:
$\frac{S_c(r,P)}{P} = \frac{p_1^2p_2^2 +
O(p_1^{\frac{3}{2}}p_2^{\frac{3}{2}})}{p_1p_2} = p_1p_2 +
O(\sqrt{p_1p_2})$
\end{enumerate}

We have proved

\begin{lem}[Conditions to Evaluate the Five Types of
Sums]\label{lemconditionstodosums} Assume the family satisfies
Equation \ref{eqexpectedAisums}. If, up to lower order terms, the
five sums (Equation \ref{eqfivesums}) are $G_P(N)
\frac{S_c(r,P)}{P}$, then the family satisfies Conditions
\ref{eqconditionsonF}.
\end{lem}

\subsection{Taylor Expansion of $G_{d,i,P}(t')$}

Fix $i$ and $d$. We calculate the first order Taylor Expansion of
$G_{d,i,P}(t') = G_P(t_i(d) + t'd^2)$. $G_{d,i,P}$ involves $t'$
only through expressions like $\frac{\log p_j}{\log C(t)}$, where
$t = t_i(d) + t'd^2$. Let $C(t) = h_m t^m + \cdots + h_0$.

The derivative of $G_{d,i,P}$ in $t'$ will involve nice functions
times factors like

\begin{eqnarray}\label{eqtaylorexpansionofgdip}
\frac{d}{dt'} \frac{\log p_j}{\log C(t)} & = & -\frac{\log
p_j}{\log^2 C(t)} \frac{d}{dt'} \log C(t_i(d) + t'd^2) \nonumber\\
& = & -\frac{\log p_j}{\log^2 C(t)} \frac{m h_m t^{m-1} d^2 +
\cdots}{h_m t^{m-1} \cdot (t_i(d) + t'd^2) + \cdots }  \nonumber\\
& \leq & \Big(\frac{10m}{|h_m|}\max_{0\leq k \leq m-1} |m-k|\cdot
|h_{m-k}|\Big) \frac{\log p_j}{\log^2 C(t)} \frac{d^2}{t_i(d) +
t'd^2}, \nonumber\\
\end{eqnarray}

provided $N$ is sufficiently large.

As $p_j \leq C(t)^\sigma$, where $\sigma$ is related to the
support of $G$, $\frac{\log p_j}{\log C(t)} \leq \sigma$. As
$C(t)$ is of size a power of $t$, we have

\begin{lem}[Taylor Expansion of
$G_{d,i,P}$]\label{lemtaylorsievea}
\begin{eqnarray}
G_{d,i,P}(t') & = & G_{d,i,P}(0) + O\Big(\frac{1}{\log N}\Big).
\end{eqnarray}
The constant above does not depend on $p_j$, $d$ or $i$.
\end{lem}

By the Mean Value Theorem $\exists \xi \in [0,t']$, corresponding
to $t_\xi = t_i(d) + \xi d^2 \in [N,2N + d^2]$ $\subset [N,2.1N]$,
such that

\begin{eqnarray}
G_{d,i,P}(t') = G_{d,i,P}(0) + \frac{d}{dt'}G_{d,i,P}\Big|_{t' =
\xi} \Big(t' - 0\Big). \end{eqnarray}

First, we have derivatives of $\frac{\log p_j}{\log C(t)}$, which
can be universally bounded from the support of $G$. Second, we
evaluate $G$ and its derivative at $2^{r_j - \kappa(r)}\frac{\log
p_j}{\log C(t_\xi)}$. We see it is sufficient to universally bound
functions like $\frac{d}{dt'}g(\frac{\log p}{\log C(t)})$.

$\log C(t_\xi) \approx \log C(N)$. Evaluating the derivative at
$\xi$, by Equation \ref{eqtaylorexpansionofgdip} we have something
bounded by $\frac{1}{\log C(t_\xi)} \frac{d^2}{t_i(d) + \xi d^2}$.
We then multiply by $t' - 0$. Thus we are bounded by
$\frac{1}{\log C(N)} \frac{t'd^2}{t_i(d) + \xi d^2}$. As $t_i(d)
\geq N$ and $t' d^2 \leq N$, the bound is at most $\frac{1}{\log
C(N)}$.

\begin{lem}[Further Taylor Expansion of $G_{d,i,P}$]\label{lemtaylorsieve}
\begin{eqnarray}
G_{d,i,P}(t') & = & G_{P}(N) + O\Big(\frac{1}{\log N}\Big).
\end{eqnarray}
The constant above does not depend on $p_j$, $d$ or $i$.
\end{lem}

The proof is similar to the previous lemma. $G_{d,i,P}(0) =
G_P\Big(t_i(d)\Big)$, $t_i(d) \in [N,N+d^2]$. Thus, to replace
$G_{d,i,P}(0)$ with $G_P(N)$ involves Taylor Expanding $G_P(t)$
around $t = N$. \hfill $\done$

This allows us to replace all the conductors of curves with $D(t)$
good with the value from $t=N$ with small error. This is very
convenient, as $G_P(N)$ has no $t'$, $i$ or $d$ dependence.
Consequently, we will be able to move it past all summations
except over primes, which will allow us to take advantage of
cancellations in $t$-sums of the $a_t(p)$'s.

\subsection{Removing the $\nu(d)||S||_\infty$ Term for $d < \log^l
N$}

\begin{eqnarray}
\sum_{t \in T(d)} S(t) = \sum_{i=1}^{\nu(d)} \sum_{t' =
0}^{[N/d^2]} S\Big(t_i(d) + t'd^2\Big) \ + \
O\Big(\nu(d)||S||_\infty \Big).
\end{eqnarray}

We show the $O\Big(\nu(d)||S||_\infty \Big)$ piece does not
contribute for $d < \log^l N$. Using Hasse to trivially bound
$||S||_\infty$ gives $2^r P^r$. We hit this with $\frac{1}{P^r}$
and sum over the primes, which will be at most $O(N^\sigma)$. We
now sum over $d < \log^l N$, getting

\begin{eqnarray}
& \ll & N^\sigma \sum_{d=1}^{\log^l N} \nu(d) \ll N^\sigma
\sum_{d=1}^{\log^l N} d^\epsilon \ll N^\sigma \log^{l(1+\epsilon)}
N. \end{eqnarray}

We then divide by the cardinality of the family, which is assumed
to be a multiple of $N$. There is no contribution for $\sigma_1 +
\sigma_2 < 1$.

\subsection{Sieving}

Let $B$ be the largest square which divides $D(t)$ for all $t$.
Recall by $t$ good we mean $D(t)$ is square-free except for primes
dividing $B$, and for $p|B$, the power of $p|D(t)$ is independent
of $t$. By Theorem \ref{thmconditionsnfbig}, possibly after
passing to a subsequence, we can approximate $t$ good by

\begin{eqnarray}
\sum_{t \in [N,2N] \atop t \ good } S(t) &=& \sum_{d=1 \atop (d,B)
= 1}^{\log^l N} \mu(d) \sum_{t \in [N,2N] \atop  D(t) \equiv 0
(d^2)} S(t) + O\Big(\sum_{t \in \mathcal{T}} S(t)\Big),
\end{eqnarray}

where the set of good $t$ is $c_{\mathcal{F}}N + o(N)$,
$c_{\mathcal{F}} > 0$, $\mathcal{T}$ is the set of $t \in [N,2N]$
such that $D(t)$ is divisible by the square of a prime $p > \log^l
N$ and $|\mathcal{T}| = o(N)$.

\subsection{Contributions from $d < \log^l N$}

We would like to use Lemma \ref{lemtaylorsieve} to replace
$G_{d,i,P}(t')$ with $G_P(N)$ plus a manageable error. This works
for pairs such as $r = (2,0)$ or $r = (2,2)$ but fails for pairs
such as $r = (1,0)$. There, we need to evaluate $\oof \sumef \oop
S(r,p)$. Replacing $\widetilde{a}_p(t)$ with $|a_t(p)| \leq
2\sqrt{p}$ gives

\begin{eqnarray}
\ll \frac{1}{|\mathcal{F}|} \frac{N}{p} \sqrt{p}, \end{eqnarray}

which is disastrous when we sum over $p$. The reason we must
trivially bound $\widetilde{a}_P(t)$ is the Taylor Expansion. We
evaluate the derivative at $\xi(t') = \xi(p_j,i,d;t')$. The
dependence of the other parameters prevents us from obtaining
complete sums (mod $P$) and using that cancellation for control.
We need to keep the cancellation from summing
$\widetilde{a}_P(t)$.

We use Partial Summation twice. Note we may always replace a
$G_{d,i,P}(t')$ with a $G_P(N)$ at a cost of $\frac{1}{\log N}$.

Let $\widetilde{A}_P(u) = \sum_{t'=0}^u \widetilde{a}_P(t')$. As
$(p_i,d) = 1$ (\textbf{this is why we are assuming $d \leq \log^l
N$ and $p_i \geq \log^l N$}), every time $t'$ increases by $P$ we
have a complete sum of the $\widetilde{a}_P$'s. Thus,

\begin{eqnarray}
\widetilde{A}_P(u) &=& \Big[\frac{u}{P}\Big]S_c(r,P) +
O\Big(P^{1+\frac{r}{2}}\Big) = \frac{u}{P}S_c(r,P) +
O\Big(P^R\Big) \nonumber\\ & & R = 1 + \frac{r}{2},\ P^R =
\prod_{j = 1 \atop r_j \neq 0}^2 p_j^{1+ \frac{r_j}{2}}.
\end{eqnarray}

In the above, the first error term is from our bound for the
incomplete sum of at most $P$ terms, each term bounded by
$\sqrt{p_1^{r_1}p_2^{r_2}} = P^{\frac{r}{2}}$. Dropping the
greatest integer brackets costs at most $S_c(r,P) = O(P^r)$. $P^r
= p_1^{r_1}p_2^{r_2}$, and $P^{1+\frac{r}{2}} = p_1^{1 +
\frac{r_1}{2}}p_2^{1 + \frac{r_2}{2}}$. As $r_j \in \{0,1,2\}$,
$r_j \leq 1 + \frac{r_j}{2}$. Thus, we may incorporate the error
from removing the greatest integer brackets into the $O(P^R)$
term.

\begin{eqnarray}\label{eqfourtermexpansion}
S(d,i,r,P) & = & \sum_{t' = 0}^{[N/d^2]} \widetilde{a}_{d,i,P}(t')
G_{d,i,P}(t') \nonumber\\ & = & \Bigg( \frac{[N/d^2]}{P} S_c(r,P)
+ O\Big(P^R \Big) \Bigg) G_{d,i,P}([N/d^2]) \nonumber\\ & & \ -
\sum_{u=0}^{[N/d^2]-1} \Bigg(\frac{u}{P}S_c(r,P) + O\Big(P^R\Big)
\Bigg)\Bigg( G_{d,i,P}(u) - G_{d,i,P}(u+1)\Bigg) \nonumber\\
S(r,P) & = & \sum_{d=1}^{\log^l N } \mu(d) \sum_{i=1}^{\nu(d)}
S(d,i,r,P) = \sum_{w=1}^4 \sum_{d=1}^{\log^l N } \mu(d)
\sum_{i=1}^{\nu(d)} S_w(d,i,r,P). \nonumber\\
\end{eqnarray}

\subsubsection{First Sum: $\frac{[N/d^2]}{P} S_c(r,P) G_{d,i,P}([N/d^2])$}

Summing over $i$ and $d$ yields

\begin{eqnarray}\label{eqtempsonerp}
S_1(r,P) &=& \sum_{d=1}^{\log^l N} \mu(d) \sum_{i=1}^{\nu(d)}
\frac{[N/d^2]}{P} S_c(r,P) G_{d,i,P}([N/d^2]) \nonumber\\ &=&
\frac{S_c(r,P)}{P} \sum_{d=1}^{\log^l N} \mu(d)
\sum_{i=1}^{\nu(d)} \Bigg[\frac{N}{d^2}\Bigg] \Bigg(G_P(N) +
O\Big(\frac{1}{\log N}\Big)\Bigg) \nonumber\\ &=&
\frac{S_c(r,P)G_P(N)}{P} \sum_{d=1}^{\log^l N} \mu(d)
\sum_{i=1}^{\nu(d)} \sum_{t'=0}^{[N/d^2]} \Bigg(1 +
O\Big(\frac{1}{\log N}\Big)\Bigg) \nonumber\\ &=&
\frac{S_c(r,P)G_P(N)}{P} \sum_{d=1}^{\log^l N} \mu(d)
\Bigg(O(\nu(d)) + \sum_{t=N \atop D(t) \equiv 0 (d^2)}^{2N} 1
\Bigg) \Bigg(1 + O\Big(\frac{1}{\log N}\Big)\Bigg) \nonumber\\ &=&
\frac{S_c(r,P)G_P(N)}{P} |\mathcal{F}| + \frac{S_c(r,P)}{P} \cdot
o(N).
\end{eqnarray}

In the last line, the error term follows from Equation
\ref{eqsizeoffamily} (which gives the $d$, $t$-sums are
$|\mathcal{F}| + o(N)$) and Lemma \ref{lemnudbound} (which gives
$\nu(d) \ll d^\epsilon$). Dividing by $|\mathcal{F}| =
c_{\mathcal{F}} N + o(N)$, the error term will not contribute when
we sum over primes, leaving us with $\frac{S_c(r,P)G_P(N)}{P}$.

\subsubsection{Second Sum: $O(P^R)G_{d,i,P}([N/d^2]) $}

Summing over $i$ and $d$ yields

\begin{eqnarray}
S_2(r,P) & \ll & \sum_{d=1}^{\log^l N} |\mu(d)|
\sum_{i=1}^{\nu(d)} P^R |G_{d,i,P}([N/d^2])| \nonumber\\ & \ll &
P^R \sum_{d=1}^{\log^l N} |\mu(d)| \sum_{i=1}^{\nu(d)}
||G||_\infty \nonumber\\ & \ll & P^R \sum_{d=1}^{\log^l N}
|\mu(d)| \sum_{i=1}^{\nu(d)} 1. \end{eqnarray}

As $\nu(d) \ll d^\epsilon$, we obtain

\begin{eqnarray}
S_2(r,P) \ \ll \ P^R \log^{l(1+\epsilon)} N \leq P^R \log^{2l} N \
= \ P^{1+\frac{r}{2}} \log^{2l} N. \end{eqnarray}

We divide by $|\mathcal{F}| = c_{\mathcal{F}}N + o(N)$, hit it
with $\frac{1}{P^r}$ and then sum over the primes. By Lemma
\ref{lemprimescontributerange}, for small support ($\sigma =
\sigma_1 + \sigma_2 < \frac{2}{3m})$ there is no contribution.

\subsubsection{Third Sum: $\sum_{u=0}^{[N/d^2]-1} \frac{u}{P}S_c(r,P)
\Big( G_{d,i,P}(u) - G_{d,i,P}(u+1)\Big) $}

We apply Partial Summation, where $a_u = G_{d,i,P}(u) -
G_{d,i,P}(u+1)$ and $b_u = \frac{u}{P}S_c(r,P)$. Thus

\begin{eqnarray}
S_3(d,i,r,P) &= & \Bigg(G_{d,i,P}(0) - G_{d,i,P}\Big([N/d^2]\Big)
\Bigg) \frac{[N/d^2]-1}{P}S_c(r,P) \nonumber\\ & & -
\sum_{u=0}^{[N/d^2]-2} \Big( G_{d,i,P}(0) - G_{d,i,P}(u+1) \Big)
\frac{1}{P}S_c(r,P). \end{eqnarray}

Using the Taylor Expansion, we gain a $\frac{1}{\log N}$ in the
first term, making it of size $\frac{S_c(r,P)}{P}
\frac{[N/d^2]}{\log N} \ll \frac{S_c(r,P)}{P}
\frac{|\mathcal{F}|}{d^2 \log N}$.

For the second term, we have $< [N/d^2]$ summands, each $\ll
\frac{1}{\log N} \frac{S_c(r,P)}{P}$. We again obtain a term of
size $\frac{S_c(r,P)}{P}\frac{|\mathcal{F}|}{d^2 \log N}$.

We sum over $i$ and $d$.

\begin{eqnarray}
S_3(r,P) &\ll & \sum_{d=1}^{\log^l N} |\mu(d)| \sum_{i=1}^{\nu(d)}
\frac{S_c(r,P)}{P}\frac{|\mathcal{F}|}{d^2 \log N} \nonumber\\
&\ll & \frac{S_c(r,P)}{P} \frac{|\mathcal{F}|}{\log N}
\sum_{d=1}^{\log^l N} \sum_{i=1}^{\nu(d)} \frac{1}{d^2}
\nonumber\\ & \ll & \frac{S_c(r,P)}{P} \frac{|\mathcal{F}|}{\log
N} \sum_{d=1}^{\log^l N} \frac{\nu(d)}{d^2}. \end{eqnarray}

As $\nu(d) \ll d^\epsilon$, $S_3(r,P) \ll \frac{S_c(r,P)}{P}
\frac{|\mathcal{F}|}{\log N}$.

\subsubsection{Fourth Sum: $\sum_{u=0}^{[N/d^2]-1} O(P^R)
\Big( G_{d,i,P}(u) - G_{d,i,P}(u+1)\Big)$}

Using the Taylor Expansion for $G_{d,i,P}(u) - G_{d,i,P}(u+1)$ is
insufficient. That gives $\frac{NP^R}{d^2\log N}$. Summing over
$i$ and $d$ is manageable, giving $O(P^R \frac{|\mathcal{F}|}{\log
N})$. Dividing by the cardinality of the family yields
$O(\frac{P^R}{\log N})$.

The problem is in summing over the primes, as we no longer have
$\frac{1}{|\mathcal{F}|}$. We multiply by $\frac{1}{P^r}$. We
recall the definitions of $r$ and $R$ and unwind the above.

Consider the case $r = (1,0)$. Then $P = p_1 = p$, $R = 1 +
\frac{r_1}{2} = \frac{3}{2}$, and $\frac{1}{P^r} = \oop$. We have

\begin{eqnarray}
\sum_{p=\log^l N}^{N^{m\sigma}} \oop \frac{p^{\frac{3}{2}}}{\log
N} \gg N^{m\sigma}. \end{eqnarray}

As $N \to \infty$, this term diverges. We need significantly
better cancellation in

\begin{eqnarray}
S_4(r,P) & = & \sum_{d=1}^{\log^l N} \mu(d) \sum_{i=1}^{\nu(d)}
\sum_{u=0}^{[N/d^2]-1} O(P^R) \Big( G_{d,i,P}(u) -
G_{d,i,P}(u+1)\Big).\nonumber\\ \end{eqnarray}

Taking absolute values and using the maximum of the $O(P^R)$ terms
gives

\begin{eqnarray}\label{eqpreboundedvariation}
S_4(r,P) &\ll & P^R \sum_{d=1}^{\log^l N} \sum_{i=1}^{\nu(d)}
\sum_{u=0}^{[N/d^2]-1} \Big| G_{d,i,P}(u) - G_{d,i,P}(u+1)\Big|.
\end{eqnarray}

The constant is independent of $P$. Taking the maximum of the
$P^R$ term involves the maximum of either the incomplete sum or
one complete sum. Using Hasse, the constant is at most
$2^{r_1+r_2}$. Thus, the constant in Equation
\ref{eqpreboundedvariation} does not depend on $P$.

If exactly one of the $r_j$'s is non-zero, then

\begin{eqnarray}
G_{d,i,P}(u) - G_{d,i,P}(u+1) = g\Big(\frac{\log p}{ \log C(t_i(d)
+ ud^2)} \Big) - g\Big(\frac{\log p}{\log C(t_i(d)
+(u+1)d^2)}\Big) \nonumber\\  \end{eqnarray}

for some Schwartz function $g$ of compact support.

If both of the $r_j$'s are non-zero, we may write $G_{d,i,P}(u)$
as the product of two functions, say $g_1$ and $g_2$. Thus

\begin{eqnarray}
G_{d,i,P}(u) = \prod_{j=1}^2 g_j\Big(\frac{\log p_j}{\log C(t_i(d)
+ ud^2)} \Big) \end{eqnarray}

Recall
\begin{eqnarray}
|a_1a_2 - b_1b_2| & = & |a_1a_2 - b_1a_2 + b_1a_2 - b_1b_2|
\nonumber\\ & \leq & |a_1a_2 - b_1a_2| + |b_1a_2 - b_1b_2| =
|a_2|\cdot |a_1 - b_1| + |b_1| \cdot |a_2 - b_2|.\nonumber\\
\end{eqnarray}

We apply the above to our function $G_{d,i,P}(u) =
g_1(d,i,p_1;u)g_2(d,i,p_2;u)$. Each $g_j(d,i,p_j;u)$ can be
bounded independently of $d$, $i$, $p_j$ and $u$, as each $g_j$ is
a Schwartz function defined in terms of the $n$-level density test
functions. Let $B = \max_j ||g_j||_\infty$. Then

\begin{eqnarray}
S_4(d,i,r,P)(u) & = & G_{d,i,P}(u) - G_{d,i,P}(u+1) \nonumber\\ &
= & \prod_{j=1 \atop r_j \neq 0}^2 g_j\Big(\frac{\log p_i}{\log
C(t_i(d) + ud^2)} \Big) - \prod_{j=1 \atop r_j \neq 0}^2
g_j\Big(\frac{\log p_j}{\log C(t_i(d) + (u+1)d^2)} \Big)
\nonumber\\ & \leq & \sum_{j=1 \atop r_j \neq 0}^2 B \cdot \Bigg|
g_j\Big(\frac{\log p_j}{\log C(t_i(d) + ud^2)}\Big) -
g_j\Big(\frac{\log p_j}{\log C(t_i(d) + (u+1)d^2)}\Big)  \Bigg|.
\nonumber\\
\end{eqnarray}

We sum the above over $u$, $i$ and $d$. Let $t_{i,d}(u) = t_i(d) +
ud^2$.

\begin{eqnarray}
S_4(r,P) &\leq & 2^rP^R\sum_{d=1}^{\log^l N} |\mu(d)|
\sum_{i=1}^{\nu(d)} \sum_{u=0}^{[N/d^2]-1} S_4(d,i,r,P)(u)
\nonumber\\ & \leq & 2^rP^R\sum_{d=1}^{\log^l N}
\sum_{i=1}^{\nu(d)} \sum_{j=1 \atop r_j \neq 0}^{2} B
\sum_{u=0}^{[N/d^2]-1} \Bigg| g_j\Big(\frac{\log p_j}{\log
C(t_{i,d}(u))}\Big) - g_j\Big(\frac{\log p_j}{\log
C(t_{i,d}(u+1))}\Big) \Bigg|.\nonumber\\ \end{eqnarray}

We show the $u$-sums are bounded independent of $p_j$, $i$, $d$,
and $N$. We may add

\begin{eqnarray}
\Bigg|g_j(0) - g_j\Big(\frac{\log p_j}{\log C(t_i(d))}\Big) \Bigg|
+ \Bigg|g_j\Big(\frac{\log p_j}{\log C(t_i(d) + [N/d^2]d^2)} \Big)
- g_j(1000\sigma)\Bigg|. \ \ \end{eqnarray}

As each $g_j$ is a Schwartz function, they are of bounded
variation. Let $x_u(d,i,p_j) = \frac{p_j}{\log N_{t_i(d)+ud^2}}$.
\textbf{As the conductors are monotone increasing}, $x_u(d,i,p_j)
> x_{u+1}(d,i,p_j)$. Thus, we have a partition of
$[0,1000\sigma]$, and we may now apply theorems on bounded
variation to bound the $u$-sum independent of $p_j$, $i$, $d$ and
$N$, obtaining $\ll 1000\sigma$.

The above is an exercise in the bounded variation of $g(x)$ on
$[0,\sigma]$. If we were to regard this as a problem in the
bounded variation of $g_{j;p_j,d,i}$ we would have $u$ ranging
over at least $\Big[0,[N/d^2]\Big]$. Even though we would gain a
$\frac{1}{\log N}$ from the derivatives, the bounded variation
bound depends on the size of the interval, which here is of length
$[N/d^2]$. We could also argue that each $g_j$ has continuous,
bounded first derivative on $[0,1000\sigma]$. By the Mean Value
Theorem, the $u$-sum is $\ll ||g_j'||_\infty \cdot |1000\sigma -
0|$.

Thus, the $u$ and the $j$-sums are universally bounded. We are
left with $\ll P^R$. Summing over $i$ and $d$ gives $\ll P^R
\log^{l(1+\epsilon)} N$. We multiply by $\frac{1}{P^r}$ and sum
over the primes. The prime sums give $N^{h(\sigma)}$; dividing by
the cardinality of the family (a multiple of $N$), we find there
is no contribution for small support. \\

\textbf{Note:} \emph{if our conductors are not monotone, we cannot
apply theorems on bounded variation. The problem is we could
transverse $[0,1000\sigma]$ (or a large subset of it) up to
$\frac{N}{d^2}$ times. This is why $S_4$ is the most difficult of
the error pieces, and why we needed to obtain polynomial
expressions for the conductors for good $t$.}

\subsubsection{Summary of Contributions for $d < \log^l N$}

\begin{lem}[Contributions for $d < \log^l N$]\label{lemdlesslogln}
Based on our Sieving Assumptions for the family (for good $D(t)$
the conductors are given by a monotone polynomial in $t$, a
positive percent of $t \in [N,2N]$ give $D(t)$ good), the main
term contribution from $d < \log^l N$ is $\frac{S_c(r,P)}{P}
G_P(N)|\mathcal{F}|$. The error terms are either of size
$\frac{S_c(r,P)}{P} o(|\mathcal{F}|)$, which won't contribute when
we sum over primes, or are such that their sum over primes will
not contribute.
\end{lem}

\subsection{Contributions from $t \in \mathcal{T}$}

\subsubsection{Preliminaries}

We are left with estimating the contributions from the troublesome
set

\begin{eqnarray}
\mathcal{T} &=& \Big\{ t \in [N,2N]: \exists d > \log^l N \
\mbox{with} \ d^2|D(t) \Big\} \end{eqnarray}

We will show in Theorem \ref{thmconditionsnfbig} that
$|\mathcal{T}| = o(N)$. By Cauchy-Schwartz

\begin{eqnarray}
\Big| \sum_{t \in \mathcal{T} } S(t) \Big|  \ \leq \ \Big(\sum_{t
\in \mathcal{T} } S^2(t) \Big)^\foh \Big(\sum_{t \in \mathcal{T} }
1 \Big)^\foh \ \leq \  \Big(\sum_{t = N}^{2N} S^2(t) \Big)^\foh
o\Big(\sqrt{N}\Big). \end{eqnarray}

We then sum over the primes, and need to show the sum over $t$ is
$O(N)$. As it stands, however, this is not sufficient to control
the error. Quick sketch: assume $S(t) = a_t(p) g(\frac{\log
p}{\log C(t)})$. Ignoring the $t$-dependence in the conductors, we
have

\begin{eqnarray}
\sum_{t = N}^{2N} S(t) &\approx & g^2\Big(\frac{\log p}{\log
C(N)}\Big)\frac{N}{p} \sum_{t (p)} a_t^2(p) \nonumber\\ & \approx
& g^2\Big(\frac{\log p}{\log C(N)}\Big) \frac{N}{p} p^2 = O(Np).
\end{eqnarray}

Taking the square-root, we hit it with $\frac{1}{p}$ and sum over
$p \leq N^\sigma$, which is not $O(\sqrt{N})$.

$S(t)$ is the product of at most two terms involving factors such
as $a_t^{r_j}(p_j)$. We hit this with factors $p_j^{-r_j}$ and sum
over $p$. Thus, instead of $S(t)$ consider $S_1(t)S_2(t)$, where
$S_j(t)$ incorporates the sum over primes to the $j^{\mbox{th}}$
power and all relevant factors.

\begin{eqnarray}
S &=& \sum_{t=N}^{2N} \Bigg[\prod_{j=1 \atop r_j \neq 0}^2
\sum_{p_j \geq \log^l N} p_j^{-r_j} g_j\Big(\frac{\log p_j}{\log
C(t)}\Big) a_t^{r_j}(p_j) \Bigg]^2 \nonumber\\ &=& \sum_{t=N}^{2N}
\prod_{w=1}^2 \prod_{j=1 \atop r_j \neq 0}^2 \sum_{p_{j_w} \geq
\log^l N} p_{j_w}^{-r_{j_w}} g_{j_w} \Big(\frac{\log p_{j_w}}{\log
C(t)}\Big) a_t^{r_{j_w}}(p_{j_w}). \end{eqnarray}

We proceed similarly as in the $d \leq \log^l N$ case, except now
there are no $d$ and $i$, and we have potentially four factors
instead of one or two. On expanding, we combine terms where we
have the same prime occurring multiple times. There are several
types of sums: four distinct primes (four factors), three distinct
primes (three factors), $\dots$, all primes the same (one factor).
We do the worst case, when there are four factors; the other cases
are handled similarly.

\subsubsection{A Specific Case: Four Distinct Primes}

Assume we have four distinct primes. Relabelling, we have
$p^{-r_i} a_t^{r_i}(p_i)$ for $i = 1$ to $4$. Let $P =
\prod_{i=1}^4 p_i$. Interchange the $t$-summation with the
$p_i$-summations. As before, we apply partial summation to
$\sum_{t=N}^{2N} \prod_{i=1}^4 a_t^{r_i}(p_i)$ $\cdot g_i(p_i,t)
p^{-r_i}$ $ = \sum_{t=N}^{2N} a(P,t) \cdot b(P,t)$, the only
change being the addition of the factors $\prod_i p^{-r_i}$. Now
$A(u) = \sum_{t=N}^u a(P,t)$ $ = \frac{u-N}{P} S_c(r,P)$ $+
O(\prod_{i=1}^4 p_i^{1+\frac{r_i}{2}})$, $S_c(r,P) = \prod_{i=1}^4
A_{r_i,\mathcal{F}}(p_i)$ by Lemma \ref{lemusefulapsum}. Let $P^R
= \prod_{i=1}^4 p_i^{1+\frac{r_i}{2}}$; the error in the partial
summation is $O(P^R)$.

As in Equation \ref{eqfourtermexpansion} we have

\begin{eqnarray}
S & = & \prod_{i=1}^4 \sum_{p_i} \sum_{t = N}^{2N} a_t^{r_i}(p_i)
\cdot p^{-r_i} G(P,t) \nonumber\\ & = & \prod_{i=1}^4 \sum_{p_i}
\Big( \frac{N}{P} S_c(r,P) + O(P^R) \Big) p_i^{-r_i} G(P,2N)
\nonumber\\ & & \ - \prod_{i=1}^4 \sum_{p_i} \sum_{u=N}^{2N-1}
\Big(\frac{u-N}{P}S_c(r,P) + O(P^R) \Big)p_i^{-r_i} \Big( G(P,u) -
G(P,u+1)\Big). \nonumber\\ \end{eqnarray}

For $r \geq 2$ by Hasse $A_{r,\mathcal{F}}(p) \leq 2^r
p^{1+\frac{r}{2}}$. For $r = 1$, $A_{1,\mathcal{F}}(p) \ll p$ by
\cite{De}. Hence $\forall r$, $A_{r,\mathcal{F}}(p) \ll p^r$.

\begin{eqnarray}
\prod_{i=1}^4 \frac{S_c(P)}{p_i} p_i^{-r_i} \ \ll \ \prod_{i=1}^4
\frac{A_{r_i,\mathcal{F}}(p_i)}{p_i^{1+r_i}} \ \ll \ \prod_{i=1}^4
\frac{p_i^{r_i}}{p_i^{1+r_i}} \ = \ \prod_{i=1}^4 \frac{1}{p_i}.
\end{eqnarray}

We can immediately handle the first sum. Inserting absolute values
yields something like

\begin{eqnarray}
\prod_{i=1}^4 \sum_{p_i} \frac{\log p_i}{\log C(2N)}
\Big|g_i\Big(\frac{\log p_i}{\log C(2N)}\Big)\Big| \frac{1}{p_i}
&\ll & \prod_{i=1}^4 O(1) \end{eqnarray}

where the last result (the sums over the primes) follows from
Corollary \ref{primeonesuma}.

Pulling out the prime factors and using partial summation again,
the third sum is handled similarly.

The second and fourth pieces are more difficult, and result in
significantly decreased support. We analyze this loss later. For
now, we need only note that the second sum is $\prod_i \sum_{p_i}
p_i^{r_i/2}$. For test functions of small support, this sum is
$o(N)$.

There is a slight obstruction in applying the same argument to the
fourth sum, namely, that $G(P,u)$ could be the product of four
factors. Similar to the identity $|a_1a_2 - b_1b_2|$ $\leq
|a_1|\cdot |a_1-b_1|$ $+ |b_1|\cdot |a_2-b_2|$, we have

\begin{eqnarray}\label{eqaiminusbiexpansion}
|a_1a_2a_3a_4 - b_1b_2b_3b_4| & \leq & \ \ \ |a_2a_3a_4|\cdot |a_1
- b_1| + |b_1a_3a_4| \cdot|a_2- b_2| \nonumber\\ & & + \
|b_1b_2a_4| \cdot|a_3 - b_3| + |b_1b_2b_3| \cdot |a_4 - b_4|
\nonumber\\ & \leq & \prod_{j=1}^4 \Big(|a_j|+|b_j|+1\Big)
\sum_{i=1}^4 |a_i - b_i|
\end{eqnarray}

The rest of the proof in this case is identical to the fourth sum
in the $d \leq \log^l N$ case.

Note: as we have always inserted absolute values before summing
over primes, it is permissible to extend from the primes are
distinct to all possible $4$-tuples.

\subsubsection{Handling the Other Cases}

The other cases (especially cases where some primes are equal) are
handled similarly. The only real change is if we have less than
four factors, and this only affects the Fourth Sum. For example,
if we have three factors instead of $4$, set $a_4 = b_4 = 1$ in
Equation \ref{eqaiminusbiexpansion}.

\subsection{Determining the Admissible Supports of the Test
Functions}

The largest errors arise from $r_i = 1$ terms, using Hasse to
trivially bound partial sums of $a_t(p)$ by $p^{3/2}$ (at most $p$
terms, each term at most $2\sqrt{p}$). Let $C(t)$ be a polynomial
of degree $m$ for good $t$. We assume all supports are at most
$\foh$ (as otherwise $p^2$ could exceed $N$, changing some of our
arguments above). In the $1$-level densities, we encounter errors
like

\begin{eqnarray}
\sum_{p=\log^l N}^{N^{\sigma m}} \frac{1}{p} \frac{\log p}{\log
N^m} g\Big( \frac{\log p}{\log N^m} \Big) p^{\frac{3}{2}} \ll
\sum_{p=\log^l N}^{N^{\sigma m}} p^\foh \ll N^{\frac{3\sigma
m}{2}}. \end{eqnarray}

We divide by $|\mathcal{F}|$, a multiple of $N$. The errors are
negligible for $\sigma < \min\Big(\frac{2}{3 m}, \foh\Big)$.

In the $2$-level density, the worst case (not including the
Cauchy-Schwartz arguments to handle the over-counting of almost
square-free numbers) was when we had two $r_i = 1$ terms. We have
two functions of support $\sigma_1$ and $\sigma_2$, and we obtain

\begin{eqnarray}
\prod_{i=1}^2 \sum_{p_i=\log^l N}^{N^{\sigma_i m}} \frac{1}{p}
\frac{\log p_i}{\log N^m} g\Big( \frac{\log p_i}{\log N^m} \Big)
p_i^{\frac{3}{2}} \ll \prod_{i=1}^2 \sum_{p_i=\log^l N}^{N^{\sigma
m}} p_i^\foh \ll N^{\frac{3(\sigma_1+\sigma_2) m}{2}}. \nonumber\\
\end{eqnarray}

We divide by a multiple of $N$ and see the errors are negligible
for $\sigma_1 + \sigma_2 < \min \Big(\frac{2}{3 m}, \foh\Big)$.
Thus, for $\sigma_1 = \sigma_2$, the support of each test function
is half that from the $1$-level density.

In applying Cauchy-Schwartz, we decrease further the allowable
support. The worst case is where we have four distinct primes with
$r_i = 1$. We sum as before, and obtain $N^{3(\sigma_1+\sigma_2)
m}$ (there is no factor of $2$ as two of the primes are associated
to test functions with support $\sigma_1$ and two to $\sigma_2$).
We take the square-root, and this must be $O(\sqrt{N})$. Thus, we
now find $\sigma_1 + \sigma_2 < \foh \frac{2}{3m}$. Setting
$\sigma_1 = \sigma_2$ yields the support is one-quarter that of
the $1$-level density.

\subsection{$1$- and $2$-Level Densities}

Assume the original family has rank $r$ over $\Qoft$. The Birch
and Swinnerton-Dyer conjecture and Silverman's Specialization
Theorem imply, for all $t$ sufficiently large, each curve's
$L$-function has $r$ family zeros at the critical point.

The Birch and Swinnerton-Dyer conjecture is only used for
interpretation purposes. The results below are derived
independently of this conjecture; however, assuming this allows us
to interpret some of the $n$-level density terms as contributions
from expected family zeros.

\begin{defi}[Non-Family Density] Let $D_{n,\mathcal{F}}^{(r)}(f)$ be the
$n$-level density from the non-family zeros (ie, the trivial
contributions from $r$ family zeros have been removed).
\end{defi}

\begin{thm}[$D_{n,\mathcal{F}}(f)$ and $D_{n,\mathcal{F}}^{(r)}(f)$,
$n = 1$ or $2$]\label{thmpreonetwodensityrank} For any
one-parameter family of rank $r$ over $\Qoft$ satisfying

\begin{enumerate}
  \item for good $t$(relative to $D(t)$), the conductors $C(t)$ are
  a monotone polynomial in $t$;
  \item up to $o(N)$, the good $t \in [N,2N]$ are obtainable by
  sieving up to $d = \log^l N$; further, the number of such $t$ is
  $|\mathcal{F}| = c_{\mathcal{F}}N + o(N)$, $c_{\mathcal{F}} >
  0$;
  \item $A_{1,\mathcal{F}}(p) = -rp + O(1)$, $A_{2,\mathcal{F}}(p) =
  p^2 + O(p^{\frac{3}{2}})$.
\end{enumerate}

Then for $f_i$ even Schwartz functions of small but non-zero
support $\sigma_i$,
\begin{eqnarray}
D_{1,\mathcal{F}}(f) &=& \hfo(0) + \foh f_1(0) + rf_1(0)
\nonumber\\ D_{1,\mathcal{F}}^{(r)}(f_1) &=& \hfo(0) + \foh f_1(0)
\end{eqnarray}
and
\begin{eqnarray}
D_{2,\mathcal{F}}(f) &=& \prod_{i=1}^2 \Bigg[ \hfi(0) + \foh
f_i(0) \Bigg] + 2\int_{-\infty}^\infty |u|\hfo(u)\hft(u)du
\nonumber\\ & &  - 2\widehat{f_1f_2}(0) - f_1(0)f_2(0) +
(f_1f_2)(0)N(\mathcal{F},-1) \nonumber\\ & & + (r^2-r)f_1(0)f_2(0)
+ r\hfo(0)f_2(0) + r f_1(0)\hft(0) \nonumber\\
D_{2,\mathcal{F}}^{(r)}(f_1) &=& \prod_{i=1}^2 \Bigg[ \hfi(0) +
\foh f_i(0) \Bigg] + 2\int_{-\infty}^\infty |u|\hfo(u)\hft(u)du
\nonumber\\ & &  - 2\widehat{f_1f_2}(0) - f_1(0)f_2(0) +
(f_1f_2)(0)N(\mathcal{F},-1).
\end{eqnarray}

Removing the contribution from $r$ family zeros, for small support
the $2$-level density of the remaining zeros agrees with
$SO(\mbox{even})$, $O$ or $SO(\mbox{odd})$ if the signs are all
even, equidistributed, or all odd. If Tate's conjecture is true
for the surface, we may interpret $r$ as the rank of $\mathcal{E}$
over $\Qoft$.

Let $m = \mbox{deg}\ C(t)$. For the $1$-level density, $\sigma <
\min(\foh,\frac{2}{3m})$. For the $2$-level density, $\sigma_1 +
\sigma_2 < \frac{1}{3m}$. For families where $\Delta(t)$ has no
irreducible factors of degree $4$ or more, the sieving is
unconditional, otherwise the results are conditional on ABC or the
Square-Free Sieve conjecture.

\end{thm}

Proof: When we sieve we obtain $\frac{S_c(r,P)G_P(N)}{P}$ plus
lower order terms. By Theorem \ref{lemconditionstodosums}, the
family satisfies Conditions \ref{eqconditionsonF}. Thus Lemma
\ref{lemrankwithassumptions} is applicable. \hfill $\done$

As remarked, we do not need to assume $A_{1,\mathcal{F}}(p) = -rp
+ O(1)$. A more cumbersome proof (using Lemma \ref{lemaepone})
handles $A_{1,\mathcal{F}}(p)$ for surfaces where Tate's
conjecture is known.

To apply Theorem \ref{thmpreonetwodensityrank}, we need

\begin{enumerate}
  \item \emph{the conductors are monotone polynomials for $D(t)$
  good;}
  \item \emph{a positive percent of $D(t)$ are good, and all but $o(N)$ of the good
  $t$ may be taken in the required arithmetic progressions;}
  \item \emph{knowledge of $A_{1,\mathcal{F}}(p)$ and $A_{2,\mathcal{F}}(p)$.}
  \end{enumerate}

For rational surfaces, by passing to a subsequence the above
conditions are satisfied. By changing $t \to ct + t_0$, Tate's
algorithm yields $C(t)$ is a monotone polynomial for $D(t)$ good
(Theorem \ref{thmcondcard}). By Theorem \ref{thmconditionsnfbig},
$|\mathcal{F}| = c_{\mathcal{F}} N + o(N)$, $c_{\mathcal{F}} > 0$
(ie, a positive percent of $D(t)$ are good). If Tate's conjecture
is true, Rosen-Silverman (Theorem \ref{thmsilvermanrosen}) gives
$A_{1,\mathcal{F}}(p)$; if $j(E_t)$ is non-constant, Michel's
Theorem (Theorem \ref{thmmichel}) gives $A_{2,\mathcal{F}}(p)$. We
have proved

\begin{thm}[Rational Surfaces Density Theorem]\label{thmonetwodensityrank}
Consider a one-parameter family of elliptic curves of rank $r$
over $\Qoft$ that is a rational surface. Assume GRH, $j(E_t)$ is
non-constant, and the ABC or Square-Free Sieve conjecture if
$\Delta(t)$ has an irreducible polynomial factor of degree at
least $4$. Let $f_i$ be an even Schwartz function of small but
non-zero support $\sigma_i$ and $m = \mbox{deg}\ C(t)$. For the
$1$-level density, $\sigma < \min (\foh,\frac{2}{3m})$. For the
$2$-level density, $\sigma_1 + \sigma_2 < \frac{1}{3m}$. Assume
the Birch and Swinnerton-Dyer conjecture for interpretation
purposes.

Let $M(t)$ be the product of the irreducible polynomials dividing
$\Delta(t)$ and not $c_4(t)$. If $M(t)$ is non-constant, then the
signs of $E_t$, $t$ good, are equidistributed as $N \to \infty$
(see \cite{He}). In this case, $N(\mathcal{F},-1) = \foh$.

After passing to a subsequence,

\begin{eqnarray}
D_{1,\mathcal{F}}(f_1) & = & \hfo(0) + \foh f_1(0) + rf_1(0) \nonumber\\
D_{1,\mathcal{F}}^{(r)}(f_1) &=& \hfo(0) + \foh f_1(0).
\end{eqnarray}
and
\begin{eqnarray}
D_{2,\mathcal{F}}(f) &=& \prod_{i=1}^2 \Bigg[ \hfi(0) + \foh
f_i(0) \Bigg] + 2\int_{-\infty}^\infty |u|\hfo(u)\hft(u)du
\nonumber\\ & &  - 2\widehat{f_1f_2}(0) - f_1(0)f_2(0) +
(f_1f_2)(0)N(\mathcal{F},-1) \nonumber\\ & & + (r^2-r)f_1(0)f_2(0)
+ r\hfo(0)f_2(0) + r f_1(0)\hft(0) \nonumber\\
D_{2,\mathcal{F}}^{(r)}(f_1) &=& \prod_{i=1}^2 \Bigg[ \hfi(0) +
\foh f_i(0) \Bigg] + 2\int_{-\infty}^\infty |u|\hfo(u)\hft(u)du
\nonumber\\ & &  - 2\widehat{f_1f_2}(0) - f_1(0)f_2(0) +
(f_1f_2)(0)N(\mathcal{F},-1).
\end{eqnarray}

The $2$-level non-family density is $SO(\mbox{even})$
($SO(\mbox{odd})$, $O$) if all curves are even (odd, the signs are
equidistributed).

Thus, for small support, the $1$- and $2$-level non-family density
agrees with the predictions of Katz and Sarnak; further, the
densities confirm that the curves' $L$-functions behave in a
manner consistent with having $r$ zeros at the critical point, as
predicted by the Birch and Swinnerton-Dyer conjecture.

\end{thm}


\section{Examples}
\setcounter{equation}{0}

\subsection{Constant Sign Families}

We consider several families where the sign of the functional
equation is always positive or negative. We verify the Katz-Sarnak
predictions, assuming only GRH.

\subsubsection{$\mathcal{F}: y^2 = x^3 + 2^4(-3)^3 (9t+1)^2$,
$9t+1$ Square-Free} \setcounter{equation}{0}

Let $\mathcal{F}: y^2 = x^3 + 2^4(-3)^3(9t+1)^2$, $t \in [N,2N]$,
$9t+1$ square-free. Note $y^2 = x^3 + 2^4(-3)^3D^2$ is equivalent
to $y^3 = x^3 + Dz^3$. Birch and Stephens \cite{BS} calculate the
sign of the functional equation for $y^3 = x^3 + Dz^3$, $D$
cube-free. It is

\begin{eqnarray} \epsilon_{E_D} & = & -w_3 \prod_{p \neq 3} w_p, \end{eqnarray}
where $w_3 = -1$ if $D \equiv \pm 1, \pm 3 (9)$ and $1$ otherwise,
$w_p = -1$ if $p|D, p \equiv 2 (3)$ and $1$ otherwise, and $D$ is
cube-free.

Consider $D = D(t) = 9t+1$. Mod $9$ it is $1$, so $-w_3$ is $1$.
Assume a prime congruent to $2$ mod $3$ divides $9t+1$. If there
were only one such prime, the remaining primes would be congruent
to $1$ mod $3$, and the product over all primes dividing $9t+1$
would be congruent to $2$ mod $3$, a contradiction. Hence the
number of primes congruent to $2$ mod $3$ dividing $9t+1$ is even.
For $9t+1$ square-free, this proves the functional equation is
even.

Applying Tate's algorithm (see \cite{Mil}), we find the conductors
$C(t)$ are $3^3(9t+1)^2$ for $9t+1$ square-free. $\delta_D = 1$,
$k = 1$, $a_k = 9$ so $\mathcal{P} = \{2, 3\}$. As $\nu(2) = 1$
and $\nu(3) = 0$, by Theorem \ref{thmconditionsnfbig}
$c_{\mathcal{F}} > 0$.

For $p \equiv 2 (3)$, $x \to x^3$ is an automorphism and $a_t(p) =
0$. Therefore in the sequel we assume all primes are congruent to
$1$ mod $3$, for any sum involving a prime congruent to $2$ mod
$3$ is zero.

For $p > 3$ and $p \equiv 1$ mod $3$, direct calculation gives

\begin{eqnarray}
A_{1,\mathcal{F}}(p) & \ = \ & 0 \nonumber\\
A_{2,\mathcal{F}}(p) & = & 2p^2 - 2p \ = \ 2p^2 + O(p).
\end{eqnarray}

From Michel's Theorem, Theorem \ref{thmmichel}, we expect
$A_{2,\mathcal{F}}(p) = p^2 + O(p^{\frac{3}{2}})$; however, his
theorem is only applicable for non-constant $j(E_t)$. As $j(E_t)$
is constant, we must directly compute $A_{2,\mathcal{F}}(p)$.
Further, as $a_t(p)$ trivially vanishes for half of the primes, we
expect and observe twice the predicted contribution at the other
primes. Finally, we will see later that the correction term to
$A_{2,\mathcal{F}}(p)$ contributes a potential lower order term to
the density functions.

By Dirichlet's Theorem for Primes in Arithmetic Progressions
(using Lemma \ref{lemprimeonesum} instead of Corollaries
\ref{primeonesuma} and \ref{primeonesumb}), we see the factors of
$2$ compensate for the restriction to primes congruent to $1$ mod
$3$, and this will be harmless in the applications.

Thus, the family satisfies the conditions of Theorem
\ref{thmonetwodensityrank} with $r = 0$. We verify (for small
support) the Katz-Sarnak predictions. As all the signs are even,
conditional only on GRH, we observe SO(even) symmetry, which is
distinguishable from SO(odd) and O symmetry.

\subsubsection{$\mathcal{F}: y^2 = x^3 \pm 4(4t+2)x$, $4t+2$
Square-Free}

Let $\mathcal{F}: y^2 = x^3 + 4(4t+2)x$, $4t+2$ square-free. We
need to study sums of $\js{x^3 \pm 4(4t+2)x}$. For $p > 2$,
changing variables by $t \to t - 2^{-1}$, $t \to \pm 16^{-1} t$,
we are led to study sums of $\js{x^3 + tx}$. If $p \equiv 3$ mod
$4$ then $\js{-1} = -1$. Changing variables $x \to -x$ shows
$a_t(p) = -\sum_{x(p)} \js{f_t(x)}$ vanishes; therefore, in the
sequel we only consider $p \equiv 1$ mod $4$.

Birch and Stephens \cite{BS} calculate the sign of the functional
equation for this family. For general $D$, $D$ not divisible by
$4$ or any fourth power, the sign of the functional equation for
the curve $y^2 = x^3 + 4Dx$ is

\begin{eqnarray}
w_\infty w_2 \prod_{p^2||D} w_p, \end{eqnarray}

where $w_\infty = $\ \mbox{sgn}$(-D)$, $w_2 = -1$ if $D \equiv
1,3,11,13$ mod $16$ and $1$ otherwise, $w_p = -1$ for $p \equiv 3
(4)$ , and $w_p = 1$ for other $p \geq 3$.

By restricting to positive, even, square-free $D$, we force the
sign of the functional equation to be odd. Hence $\epsilon_D = -1$
if $D = 4t+2$, $D$ square-free. If we had taken $D = -(4t+2)$,
$4t+2$ square-free, we would have found $\epsilon_D = +1$.

From Tate's algorithm, for $D(t) = \pm(4t+2)$ square-free, $C(t) =
2^6(4t+2)^2$. $\delta_D = 1$, $k = 1$, $a_k = 4$ so $\mathcal{P} =
\{2\}$. As $\nu(2) = 0$, by Theorem \ref{thmconditionsnfbig}
$c_{\mathcal{F}} > 0$.

For $p > 2$ and $p \equiv 1$ mod $4$, direct calculation gives

\begin{eqnarray}
A_{1,\mathcal{F}}(p) & \ = \ & 0 \nonumber\\
A_{2,\mathcal{F}}(p) &=& 2p^2 - 2p \ = \ 2p^2 + O(p).
\end{eqnarray}

For the family $\mathcal{F}_{\pm}: y^2 = x^3 \pm 4(4t+2)x$, $4t+2$
square-free, all curves in $\mathcal{F}_-$ have even sign, all
curves in $\mathcal{F}_+$ have odd sign. The families satisfy the
conditions of Theorem \ref{thmonetwodensityrank} with $r = 0$. We
verify (for small support) the Katz-Sarnak predictions. As all the
signs are even (odd), conditional only on GRH, we observe SO(even)
(SO(odd)) symmetry.

\subsubsection{$\mathcal{F}: y^2 = x^3 + tx^2 - (t+3)x + 1$}

For this family (due to Washington)

\begin{eqnarray}
c_4(t) &=& 2^4(t^2 + 3t + 9) \nonumber\\ \Delta(t) &=& 2^4(t^2 +
3t + 9)^2 \nonumber\\ j(E_t) &=& 2^8(t^2 + 3t + 9). \end{eqnarray}

Washington (\cite{Wa}) proved the rank is odd for $t^2 + 3t + 9$
square-free, assuming the finiteness of the Tate-Shafarevich
group. Rizzo \cite{Ri} proved the rank is odd for all $t$. While
$j(E_t)$ is non-constant, $M(t) = 1$ ($M(t)$ is the product of all
irreducible polynomials dividing $\Delta(t)$ but not $c_4(t)$).
Thus, Helfgott's results on equidistribution of sign are not
applicable.

For sieving convenience, we replace $t$ with $12t+1$. Let $D(t) =
144t^2 + 60t + 13$. Tate's algorithm yields for $D(t)$
square-free, $C(t) = 2^3(144t^2 + 60t + 13)$.

$\delta_D = -2^4 3^5$, $k = 2$, $a_k = 2^4 3^2$ so $\mathcal{P} =
\{2,3\}$. $D(t)$ is a primitive integral polynomial. For $p\notdiv
6$ the number of incongruent solutions of $D(t) \equiv 0$ mod
$p^2$ equals the number of incongruent solutions of $D(t) \equiv
0$ mod $p$ (see \cite{Nag}). As $\nu(2) = \nu(3) = 0$, by Theorem
\ref{thmconditionsnfbig}, $c_{\mathcal{F}} > 0$.

Direct calculation gives

\bea A_{1,\mathcal{F}}(p) & \ = \ & -p\Bigg[1 + \js{-1}\Bigg].
\end{eqnarray}

Hence $A_{1,\mathcal{F}}(p)$ is $-2p$ for $p \equiv 1 (4)$ and $0$
for $p \equiv 3 (4)$. By Theorem \ref{thmsilvermanrosen}, the rank
over $\Qoft$ is $1$.

As $j(E_t)$ is non-constant, by Michel's Theorem
$A_{2,\mathcal{F}}(p) = p^2 + O(p^{\frac{3}{2}})$.

The conditions of Theorem \ref{thmonetwodensityrank} are satisfied
with $r = 1$. We again verify the Katz-Sarnak predictions: there
are two pieces to our densities. The first equals the contribution
from $1$ zero at the critical point; the second agrees with
SO(odd) for small support.

\subsection{Rational Families}

We give two examples of rational families of elliptic curves over
$\Q(t)$. See \cite{Mil} for proofs, as well as a new method to
generate rational families of moderate rank.

\subsubsection{Rank $1$ Example}

Consider the rational family $y^2 = x^3 + 1 + tx^2$.

\begin{eqnarray}
c_4(t) &=& 16t^2 \nonumber\\ \Delta(t) &=& -16(4t^3 + 27)
\nonumber\\ j(E_t) &=& -256 \frac{t^6}{4t^3+27} \nonumber\\ M(t)
&=& 4t^3 + 27. \end{eqnarray}

If we replace $t$ with $6t+1$, we can easily calculate the
conductors for $D(t) = 4(6t+1)^3 + 27$ square-free. In \cite{Mil}
we show $C(t) = 2^2\Big(4(6t+1)^3 + 27\Big)$ for $D(t)$
square-free. By Hooley (\cite{Ho}, Theorem $3$, page $69$), as
$D(t)$ is an irreducible polynomial of degree $3$,
$c_{\mathcal{F}} > 0$.

Direct calculations \cite{Mil} gives $A_{1,\mathcal{F}}(p) = -p$,
and a more involved calculation gives $A_{2,\mathcal{F}}(p)$ $ =
p^2 - 3p h_{3,p}(2)$ $- 1 + p\sum_{x (p)} \js{4x^3+1}$  $= p^2 +
O(p^{\frac{3}{2}})$, where $h_{3,p}(2)$ is one if $2$ is a cube
mod $p$ and zero otherwise. Note this shows Michel's bound for
$A_{2,\mathcal{F}}(p)$ is sharp.

As $j(E_t)$ and $M(t)$ are non-constant, we expect the signs to be
equidistributed.

The Rational Surfaces Density Theorem is applicable, and we obtain
orthogonal symmetry for the density of the non-family zeros.

\subsubsection{Rank $6$ Example}

We give a more exotic example. See \cite{Mil} for the details. Let

\begin{equation}
\begin{array}{ccrr}
A & = & 8916100448256000000 & \nonumber\\ B & = &
-811365140824616222208 & \nonumber\\ C & = &
26497490347321493520384 & \nonumber\\ D & = &
-343107594345448813363200 & \nonumber\\ a & = & 16660111104 &
\nonumber\\ b & = & -1603174809600 & \nonumber\\ c & = &
2149908480000 &
\end{array} \end{equation}

The rational family $y^2 = x^3 t^2 + 2g(x)t - h(x)$, $g(x) = x^3 +
ax^2 + bx + c$ and $h(x) = (A-1)x^3 + Bx^2 + Cx + D$, has
$A_{1,\mathcal{F}}(p) = -6p + O(1)$ for $p$ large. Therefore, the
family has rank $6$ over $\Qoft$. Writing in Weierstrass normal
form yields

\begin{eqnarray}
y^2 & = & x^3 + (2at-B)x^2 + (2bt-C)(t^2 + 2t - A + 1)x
\nonumber\\ & & \ \ \ \ \ \ + (2ct - D)(t^2 + 2t - A + 1)^2
\nonumber\\ c_4(t) &=& 2^{19}3^7 7^1 13^1 (1475 t^3 + \cdots -
\ 7735999878503076170786750620939) \nonumber\\
c_6(t) & = & -2^{25}3^{11}(625 t^5 + \cdots ) \nonumber\\ j(E_t) &
= & \frac{50141357421875 t^{9} + \cdots }{-1171875t^{10} + \cdots } \nonumber\\
\Delta(t) & = & -2^{44}3^{18} 5^6(75 t^{10} + \cdots).
\end{eqnarray}

This is a rational surface, $j(E_t)$ and $M(t)$ are non-constant.
Thus, by the Rational Surfaces Density Theorem, we verify the
Katz-Sarnak predictions for a family of rank $6$ over $\Q(t)$!


\section{Summary and Future Work}
\setcounter{equation}{0}

Our main result is that, modulo standard conjectures, the
fluctuations of the non-family low lying zeros in one-parameter
families of elliptic curves agree with the Katz-Sarnak
conjectures. Further, a family of rank $r$ over $\Q(t)$ has a
density correction which equals the contribution of $r$ zeros at
the critical point, providing further evidence for the Birch and
Swinnerton-Dyer conjecture.

We have found four families where the observed density agrees with
the density of one (and only one) symmetry group. As expected, the
first piece equals the contribution from $r$ zeros at the critical
point (where $r$ is the geometric rank of the family), and the
second equals SO(even) if all curves have even sign and SO(odd) if
all curves have odd sign.

For these four families, we assumed only GRH. We are able to
unconditionally handle the dependence of the conductors on $t$,
the signs of the functional equations, and the error terms.

In general, the greatest difficulty is handling the variation in
the conductors. Unlike other families investigated (\cite{ILS},
\cite{Ru}), the conductors of elliptic curves vary wildly in a
given family. If the discriminant $\Delta(t)$ has an irreducible
factor of degree $4$ or greater, either ABC or the Square-Free
Sieve Conjecture must be assumed to perform the necessary sieving;
if all irreducible factors are of degree at most $3$, the sieving
is unconditional.

The crucial observation is that, if we sieve to a positive percent
subset where the conductors are monotone, then we can bound the
error terms. Note the extreme delicacy of our arguments: for
conductors of size $\log N$, we cannot bound the error terms if
the conductors range from $\log N - \log c$ to $\log N + \log c$
for some constant $c$.

It was observed in \cite{Mil} that in every family where
$A_{2,\mathcal{F}}(p)$ can be directly calculated,

\begin{equation}
A_{2,\mathcal{F}}(p) = p^2 + h(p) - m_{\mathcal{F}} p + O(1),
\end{equation}

where $h(p)$ is of size $p^{\frac{3}{2}}$ and averages to zero,
and $m_{\mathcal{F}}$ is a positive constant, often different for
different families.

We have shown all rational families (with the same distribution of
signs) have equal $1$ and $2$-level densities. We can, however,
try to expand the densities in powers of $\frac{1}{\log N}$. The
different $m_{\mathcal{F}} p$ terms will lead to potential
corrections to the densities of size $\frac{1}{\log N}$, giving
the exciting possibility of distinguishing different families by
lower order corrections to the common densities.

Unfortunately, the size of the errors in the $1$ and $2$-level
densities are $O\Big(\frac{\log \log N}{\log N}\Big)$; thus, a
significantly more delicate analysis is needed before we can
expand the densities.

\appendix   


\section{Sieving Families of Elliptic Curves}
\setcounter{equation}{0}

Given a one-parameter family of elliptic curves $E_t$, we need to
control the conductors $C(t)$ to determine the $1$- and $2$-level
densities. Let the curves have discriminants $\Delta(t)$, and let
$D(t)$ be the product of the irreducible polynomial factors of
$\Delta(t)$.

$D(t)$ may always be divisible a fixed square; let $B$ be the
largest square dividing $D(t)$ for all $t$. We prove in Theorem
\ref{thmcondcard} that for a rational elliptic surface, by passing
to a subsequence $\tau = c_1t + c_0$, for $\frac{D(\tau)}{B}$
square-free, $C(t)$ is given by a polynomial in $\tau$. Call such
$t$ (or $D(t)$ or $\tau$) good.

In order to evaluate the sums of $\prod_i a_t^{r_i}(p_i)$, it is
necessary to restrict $t$ to arithmetic progressions; however,
restricting to $t$ good ($\frac{D(\tau)}{B}$ square-free) does not
yield $t$ in arithmetic progressions.

We overcome this difficulty by doing a partial sieve with good
bounds on over-counting. For notational convenience, we consider
the case where $B = 1$ below, and indicate how to modify for
general $B$.

Let $S(t)$ be some quantity associated to our family which we
desire to sum over $\mathcal{T}_{sqfree}$, where

\begin{eqnarray}
\mathcal{T}_{sqfree} &=& \Big\{t \in [N,2N]:\ D(t) \ \mbox{is
sqfree} \Big\} \nonumber\\ \mathcal{T}_N &=& \Big\{t \in [N,2N]:\
d^2 \notdiv D(t) \ \mbox{for} \ 2 \leq d \leq \log^l N \Big\}.
\end{eqnarray}

Clearly $\mathcal{T}_{sqfree} \subset \mathcal{T}_N$. We show
$\mathcal{T}_N$ is a union of arithmetic progressions, and
$|\mathcal{T}_N - \mathcal{T}_{sqfree}| = o(N)$.

The main obstruction is estimating the number of $t \in [N,2N]$
with $D(t)$ divisible by the square of a prime $p \geq \log^l N$.
If $k =\mbox{deg}\ D(t)$,

\begin{eqnarray}
\sum_{D(t) \ sqfree \atop t \in [N,2N]} S(t) &=&
\sum_{d=1}^{N^{k/2}} \mu(d) \sum_{D(t) \equiv 0 (d^2) \atop t \in
[N,2N] } S(t) \nonumber\\ &=& \sum_{d=1}^{\log^l N} \mu(d)
\sum_{D(t) \equiv 0 (d^2) \atop t \in [N,2N] } S(t) + \sum_{d \geq
\log^l N}^{N^{k/2}} \mu(d) \sum_{D(t) \equiv 0 (d^2) \atop t \in
[N,2N] } S(t).\nonumber\\ \end{eqnarray}

For $k > 3$, the second piece is too difficult to estimate --
there are too many $d$ terms ($d$ runs to $N^{k/2}$). If all the
irreducible factors of $D(t)$ are of degree at most $3$, the
second piece is small. For factors of degree at most $2$, this
follows immediately, while for factors of degree $3$ it follows
from Hooley (\cite{Ho}). For larger degrees, we need the ABC
conjecture (or one of its consequences, the Square-Free Sieve
conjecture).

\subsection{Incongruent Solutions of Polynomials}

Recall the following basic facts (see, for example, \cite{Nag})
for an integral polynomial $D(t)$ of degree $k$ and discriminant
$\delta$:

\begin{enumerate}
  \item Let $p$ be a prime not dividing the coefficient of
$x^k$. Then $D(t) \equiv 0$ mod $p$ has at most $k$ incongruent
solutions.
  \item Let $D(t) \equiv 0$ mod $p_i^{\alpha_i}$ have $\nu_i$ incongruent
solutions. If the primes are distinct, there are $\prod_{i=1}^r
\nu_i$ incongruent solutions of $D(t) \equiv 0$ mod $\prod_{i=1}^r
p_i^{\alpha_i}$.
  \item Suppose $p \notdiv \delta$. Then the number of incongruent solutions
  of $D(t) \equiv 0$ mod $p$ equals the number of incongruent solutions of
$D(t) \equiv 0$ mod $p^{\alpha}$.
\end{enumerate}

\begin{defi} Let $\nu(d)$ be the number of incongruent solutions of $D(t)
\equiv 0$ mod $d^2$.
\end{defi}

\begin{lem}\label{lemnudbound} For $d$ square-free, $\nu(d) \ll d^\epsilon$.
\end{lem}

The proof combines the above facts with the standard bound of the
divisor function, $\tau(d) \ll d^{\epsilon}$.

\subsection{Common Prime Divisors of Polynomials}

\begin{lem}\label{lemftgt} Let $f(t)$ and $g(t)$ be integer
polynomials with no non-constant factors over $\Z[t]$. Then
$\exists c$ (independent of $t$) such that if $p$ divides both
$f(t)$ and $g(t)$, then $p|c$. In particular, $f(t)$ and $g(t)$
have no common large prime divisors.
\end{lem}

Proof: Euclid's algorithm.

\subsection{Calculating $|\mathcal{T}_N|$}

\begin{eqnarray}
\sum_{t \in \mathcal{T}_N} 1 &=& \sum_{d=1}^{\log^l N} \mu(d)
\sum_{D(t) \equiv 0 (d^2) \atop t \in [N,2N] } 1. \end{eqnarray}

There are $\frac{N}{d^2}\nu(d) + O(\nu(d))$ solutions to $D(t)
\equiv 0$ mod $d^2$ for $t \in [N,2N]$. By Lemma
\ref{lemnudbound}, $\nu(d) \ll d^\epsilon$ for square-free $d$.
Thus

\begin{eqnarray}
|\mathcal{T}_N| &=& \sum_{d=1}^{\log^l N} \mu(d)
\Bigg[\frac{N}{d^2}\nu(d) + O(\nu(d)) \Bigg] = N
\sum_{d=1}^{\log^l N} \frac{\mu(d) \nu(d)}{d^2} +
O(\log^{l(1+\epsilon)} N).\nonumber\\ \end{eqnarray}

As $\nu(d) \ll d^\epsilon$ for square-free $d$,

\begin{eqnarray}
\Big| \prod_{p < \log^l N} \Big( 1 - \frac{\nu(p)}{p^2}\Big) -
\sum_{d=1}^{\log^l N} \frac{\mu(d)\nu(d)}{d^2} \Big| & \ll &
\sum_{d=\log^l N}^{\infty} \frac{d^\epsilon}{d^2} \ll
\frac{1}{\log^{l(1-\epsilon)} N}. \nonumber\\ \end{eqnarray}

Therefore

\begin{eqnarray}
|\mathcal{T}_N| &=& N \prod_{p < \log^l N} \Big( 1 -
\frac{\nu(p)}{p^2}\Big) + O\Big( \frac{N}{\log^{l(1-\epsilon)} N}
\Big) + O(\log^{l(1-\epsilon)} N). \nonumber\\ \end{eqnarray}

We may take the product over all primes with negligible cost as

\begin{eqnarray}
1 - \prod_{p \geq \log^l N} \Big( 1 - \frac{\nu(p)}{p^2} \Big) \ll
\sum_{n \geq \log^l N} \frac{n^\epsilon}{n^2} \ll
\frac{1}{\log^{l(1-\epsilon)} N}. \end{eqnarray}

We have shown

\begin{lem} $\mathcal{T}_N = \{t \in
[N,2N]: d^2 \notdiv D(t)$ for $2 \leq d \leq \log^l N \}$.
\begin{eqnarray}\label{eqTN}
|\mathcal{T}_N| &=& N \prod_p \Big( 1 - \frac{\nu(p)}{p^2}\Big) +
O\Big( \frac{N}{\log^{l(1-\epsilon)} N} \Big).
\end{eqnarray}
\end{lem}

\subsection{Estimating $\mathcal{T}_{sqfree}$}

Assuming the ABC conjecture, Granville (\cite{Gr}, Theorem 1)
proves the number of $t \in [N,2N]$ such that $D(t)$ is
square-free is

\begin{eqnarray}\label{eqtsqfree}
|\mathcal{T}_{sqfree}| &=& N \prod_p \Big( 1 - \frac{\nu(p)}{p^2}
\Big) + o(N).
\end{eqnarray}

Again, if the degree of $D(t)$ is at most $3$, the ABC conjecture
is not needed. The family has a positive percent of $t$ giving
$D(t)$ square-free (as we are assuming no square divides $D(t)$
for all $t$, no $\nu(p) = p^2$, hence the product can be bounded
away from $0$).

\subsection{Evaluation of $|\mathcal{T}_N - \mathcal{T}_{sqfree}|$ and Applications}

By Equations \ref{eqTN} and \ref{eqtsqfree}, as
$\mathcal{T}_{sqfree} \subset \mathcal{T}_N$, we have $
|\mathcal{T}_N - \mathcal{T}_{sqfree}| = o(N)$.

We have proved

\begin{eqnarray}
\sum_{t \in [N,2N] \atop D(t) \ sqfree} S(t) &=& \sum_{t \in
\mathcal{T}_N } S(t) + O\Big( \sum_{t \in \mathcal{T}} S(t) \Big)
\nonumber\\ &=& \sum_{d=1}^{\log^l N} \mu(d) \sum_{D(t) \equiv 0
(d^2) \atop t \in [N,2N] } S(t) + O\Big( \sum_{t \in \mathcal{T}}
S(t) \Big).
\end{eqnarray}

We use arithmetic progressions to handle the piece with $d \leq
\log^l N$, and Cauchy-Schwartz to handle $t \in \mathcal{T}$.

\begin{eqnarray}
\sum_{t \in \mathcal{T}} S(t) \ \ll \  \Big( \sum_{t \in
\mathcal{T}} S^2(t) \Big)^{\foh} \Big( \sum_{t \in \mathcal{T}} 1
\Big)^{\foh} \ \ll \  \Big( \sum_{t \in [N,2N]} S^2(t)
\Big)^{\foh} o\Big(\sqrt{N}\Big).\nonumber\\
\end{eqnarray}

If we can show $\sum_{t=N}^{2N} S^2(t) = O(N)$, then the error
term is negligible as $N \to \infty$.

\subsection{Conditions Implying $|\mathcal{F}| = c_{\mathcal{F}}
N + o(N)$, $c_{\mathcal{F}} > 0$}

Assume no square divides $D(t)$ for all $t$. The number of $t \in
[N,2N]$ with $D(t)$ not divisible by $d^2$, $d \leq \log^l N$, is
$N\prod_p \Big(1 - \frac{\nu(p)}{p^2}\Big) + o(N)$. Let $D(t) =
\prod_i D_i^{r_i}(t)$, $D_i(t)$ irreducible. By multiple
applications of Lemma \ref{lemftgt}, $\exists c$ such that
$\forall t$, there is no prime $p > c$ which divides two of the
$D_i(t)$. Thus, if $D(t)$ is divisible by $p^2$ for a large prime,
one of the factors is divisible by $p^2$. As there are finitely
many factors, it is sufficient to bound by $o(N)$ the number of $t
\in [N,2N]$ with $p^2 |D(t)$ for a large prime for irreducible
$D(t)$.

Let $|\mathcal{F}|$ equal the number of $t \in [N,2N]$ with $D(t)$
square-free. Let $c_{\mathcal{F}} = \prod_{p \leq \log^l N} \Big(1
- \frac{\nu(p)}{p^2}\Big)$. We have seen extending the product to
all primes costs $O(\frac{1}{\log^{l(1-\epsilon)} N})$. Thus, we
need only bound $c_{\mathcal{F}}$ away from zero.

Let $D(t) = a_k t^k + \cdots + a_0$ with discriminant $\delta$.
For $p \notdiv a_k \delta$, $\nu(p) \leq k$.

Let $\mathcal{P}$ be the set of primes dividing $a_k \delta$ and
all primes at most $\sqrt{k}$. The contribution from $p \not\in
\mathcal{P}$ is bounded away from $0$. Therefore, if $\nu(p) <
p^2$ for $p|a_k \delta$ and $p \leq \sqrt{k}$, then
$c_{\mathcal{F}} > 0$.

If $D(t)$ is divisible by a square for all $t$, the above
arguments fail. Let $P$ be the largest product of primes such that
$\forall t$, $P^2 | D(t)$. By changing variables $\tau \to P^m t +
t_0$, for $m$ sufficiently large, $D(\tau)$ is divisible by fixed
powers of $p|P$, depending only on $D(t_0)$. Thus, instead of
sieving to $D(t)$ square-free, we sieve to $D(\tau)$ square-free
except for primes dividing $P$.

Let $\delta_\tau$ denote the new discriminant. As the discriminant
is a product over the differences of the roots, $t_0$ does not
change the discriminant, and $P^m$ rescales by a power of $P$.
Thus, $\delta_\tau = P^M \delta$. Further, the new leading
coefficient is $P^{mk} a_k$. Thus, for $p \notdiv P$, our previous
arguments are still applicable, except we are no longer sieving
over $p | P$. We have shown

\begin{thm}[Conditions on $D(t)$ implying $|\mathcal{F}| =
c_{\mathcal{F}} N + o(N)$]\label{thmconditionsnfbig} Assume no
square divides $D(t)$ for all $t$. Let $\mathcal{P}$ be the set of
primes dividing $a_k \delta$ and all primes at most $\sqrt{k}$. If
$\forall p \in \mathcal{P}$, $\nu(p) \leq p^2 - 1$, then
$|\mathcal{F}| = c_{\mathcal{F}} N + o(N)$, $c_{\mathcal{F}} > 0$.
If $\forall t$, $B^2 | D(t)$ ($\exists p \in \mathcal{P}$, $\nu(p)
= p^2$), let $P$ be the product of all primes either in
$\mathcal{P}$ or dividing $B$. By changing variables to $\tau =
P^m t + t_0$ for $m$ large and sieving to $D(\tau)$ square-free
except for $p|P$ (where $\forall t$, the power of $p|P$ dividing
$D(t)$ is constant), we again obtain $|\mathcal{F}| =
c_{\mathcal{F}} N + o(N)$, $c_{\mathcal{F}} > 0$. In this case,
$c_{\mathcal{F}}$ no longer includes factors from $p|P$.

If all irreducible factors of $D(t)$ have degree at most $3$,
these results are unconditional; if there is an irreducible factor
with degree at least $4$ these results are conditional, and a
consequence of the ABC or Square-Free Sieve conjecture.

Further, let $\mathcal{T} = \{t \in [N,2N]: \exists d > \log^l N \
\mbox{with}\ d^2|D(t)\}$. Then $\mathcal{T} = o(N)$.

\end{thm}


\section{Handling the Conductors $C(t)$}
\setcounter{equation}{0}

For many families of elliptic curves, by sieving to a positive
percent subsequence of $t$ we obtain a sub-family where the
conductors are a monotone polynomial in $t$. In particular, we
prove this for all rational surfaces.

Tate's algorithm (see \cite{Cr}, pages $49-52$) allows us to
calculate the conductor $C(t)$ for an elliptic curve $E_t$ over
$\Q$:

\begin{eqnarray}
C(t) = \prod_{p|\Delta(t)} \ p^{f_p(t)},
\end{eqnarray}

where for $p > 3$, if the curve is minimal for $p$ then $f_p(t) =
0$ if $p\notdiv \Delta(t)$, $1$ if $p|\Delta(t)$ and $p\notdiv
c_4(t)$, and $2$ if $p|\Delta(t)$ and $p|c_4(t)$. If $p > 3$ and
$p^{12} \notdiv \Delta(t)$, then the equation is minimal at $p$.
See \cite{Si1}.

Let $\Delta(t) = d\Delta_1(t) \Delta_2(t)$, where
$\Big(\Delta_2(t), c_4(t)\Big) = 1$ and $\Delta_1(t)$ is the
product of powers of irreducible polynomials dividing $\Delta(t)$
and $c_4(t)$. By possibly changing $d$, we may take $\Delta_i(t)$
primitive. Let $D_i(t)$ be the product of all irreducible
polynomials dividing $\Delta_i(t)$, $D(t) = D_1(t)D_2(t)$.

For $t$ with $D(t)$ square-free except for small primes, $C(t) =
D_1^2(t) D_2(t)$ if $\Delta(t)$ has no irreducible polynomial
factor occurring at least $12$ times (except for corrections from
the small primes). Hence, while $f_p(t)$ may vary, the product of
$p^{f_p(t)}$, except for a finite set of primes, is well behaved.

Let

\begin{eqnarray}
\mathcal{P}_0 &=& \{p: p \leq \mbox{deg}\ \Delta(t) \} \ \cup \
\{p: p| cd\}, \ \ \ P_0 = \prod_{p \in \mathcal{P}_0} p.
\end{eqnarray}

The idea is that while for such $p$, $f_p(t)$ may vary, by
changing variables from $t$ to $P_0^m t + t_1$ for some enormous
$m$, for $p \in \mathcal{P}_0$, $f_p(P_0^m t + t_1) = f_p(t_1)$.
Thus, for this subsequence and these primes, $f_p(t)$ is constant.

We need two preliminary results. First, given a finite set of
primes $\mathcal{P}_0$, we may find an $m$ and a $t_1$ such that
for those primes, $f_p(P_0^m t + t_1)$ is constant. Second, Lemma
\ref{lemftgt}: given two polynomials with no non-constant factors
over $\Q$, there is a finite set of primes $\mathcal{P}_2$ such
that if $\exists t$ such that $\exists p$ dividing both
polynomials, then $p \in \mathcal{P}_2$.

\subsection{$f_p(t)$, $p \in \mathcal{P}_0$}

Consider the original family of elliptic curves

\begin{eqnarray}
E_t: \ y^2 + a_1(t)xy + a_3(t)y = x^3 + a_2(t)x^2 + a_4(t)x +
a_6(t). \end{eqnarray}

Assume $\Delta(t)$ is not identically zero. Choose $t_1$ such that
$\forall t \ge t_1$, $\Delta(t) \neq 0$. Apply Tate's algorithm to
$E_{t_1}$. If the initial equation was non-minimal for $p$, we
change coordinates by $T(0,0,0,p)$ (see \cite{Cr}) and restart.
After finitely many passes, Tate's algorithm terminates.

In determining $f_p(t_1)$, assume we passed through Tate's
algorithm $L_{t_1}(p)$ times. For each prime $p$, after possibly
many coordinate changes, one of the following conditions held:
$p\notdiv \Delta$, $p \notdiv c_4$, $p^2\notdiv a_6$, $p^3 \notdiv
b_8$, $p^3 \notdiv b_6$, $p\notdiv w(a_2,a_4,a_6)$, $p\notdiv
xa_3^2(a_3) + 4xa_6(a_6)$, $p\notdiv xa_4^2(a_4) -
4xa_2(a_2)xa_6(a_6)$, $p^4 \notdiv a_4$, $p^6 \notdiv a_6$, and
every function is polynomial in the $a_i$'s. Thus, after possibly
many coordinate changes, some polynomial (with integer
coefficients) of the $a_i$'s is not divisible by either $p$,
$p^2$, $p^3$, $p^4$ or $p^6$.

Consider $\tau = P_0^m t + t_1$. For $m$ enormous, $f_p(\tau) =
f_p(t_1)$ for $p \in \mathcal{P}_0$ because in Tate's algorithm,
we only need the values modulo a power of $p$. We have

\begin{eqnarray}
a_i(\tau) = a_i(P_0^m t + t_1)  = P_0^m t\widehat{a}_i(P_0^mt) +
a_i(t_1) = \widetilde{a}_i(t) + a_i(t_1). \end{eqnarray}

If $m$ is sufficiently large, we can ignore $\widetilde{a}_i(t)$
in all equivalence checks, as for these powers of $p$,
$\widetilde{a}_i(t) \equiv 0$. Let

\begin{eqnarray}
n_t(p) &=& \mbox{ord}\Big(p,\Delta(t)\Big) \nonumber\\ n &=&
\max_{p \in \mathcal{P}_0}  n_{t_1}(p)  \nonumber\\ L &=& \max_{p
\in \mathcal{P}_0}  L_{t_1}(p). \end{eqnarray}

We prove $f_p(\tau) = f_p({t_1})$ for large $m$. How large must
$m$ be? Excluding lines $42 - 65$, on each pass through Tate's
algorithm we sometimes divide our coefficients by powers of $p$:
up to $p^2$ on lines $26$ and $30$, up to $p^3$ on line $34$, up
to $p^4$ on line $69$, and $p^{12}$ on line $80$. Over-estimating,
we divide by at most $p^{2\cdot2 + 1\cdot 3 + 1\cdot 4 + 1\cdot
12}$ $= p^{23}$.

For lines $42-65$, we have a loop which can be executed at most $n
+ 4$ times. We constantly divide by increasing powers of $p$; the
largest power is the last time through the loop, which is at most
$p^{2(n+6)}$. As we pass through this loop at most $n+4$ times, we
divide by at most $p^{2n^2 + 20n + 48}$.

Thus, on each pass we have divisions by at most $p^{2n^2 + 20n +
48 + 23}$. As we loop through the main part of Tate's algorithm at
most $L$ times, we have divisions by at most $p^{(2n^2+20n+71)L}$.
If $m > (2n^2+20n+71)L$, then $\forall t$, none of the
$\widetilde{a}_i(t) = P_0^m t \widehat{a}_i(t)$ terms affect any
congruence. Significantly smaller choices of $m$ work: many of the
divisions (for example, from lines $42 - 65$) arise only once.

\subsection{Rational Surfaces I}

\subsubsection{Preliminaries}

Recall an elliptic surface $y^2 = x^3 + A(t)x + B(t)$ is rational
iff one of the following is true: $(1)$ $0 < \max\{3 \mbox{deg}\
A(t), 2\mbox{deg}\ B(t)\} < 12;$ $(2)$ $3\mbox{deg}\ A(t) =
2\mbox{deg}\ B(t) = 12$ and $\mbox{ord}_{t=0}t^{12} \Delta(t^{-1})
= 0$. See \cite{RSi}, pages $46-47$ for more details.

Assume we are in case $(1)$. No non-constant polynomial of degree
$11$ or more divides $\Delta(t)$; however, a twelfth or higher
power of a prime might divide $\Delta(t)$. Let $k = \mbox{deg}\
\Delta(t)$, and write

\begin{eqnarray}
\Delta(t) &=& d  \Delta_1(t) \Delta_2(t) \nonumber\\ c_4(t) &=& c
\gamma_1(t) \gamma_2(t) \nonumber\\ \mathcal{P}_0 &=& \{p: p \leq
\mbox{deg}\ \Delta(t) \} \ \cup \ \{p: p| cd\}, \ \ \ P_0 =
\prod_{p \in \mathcal{P}_0} p. \end{eqnarray}

where $\Delta_1(t)$ through $\gamma_2(t)$ are primitive
polynomials, $\Delta_1(t)$ and $\gamma_1(t)$ are divisible by the
same non-constant irreducible polynomials, and $\Delta_2(t)$ and
$c_4(t)$ are not both divisible by any non-constant polynomial.

Let $D_i(t)$ be the product of all non-constant irreducible
polynomials dividing $\Delta_i(t)$, and similarly for $c_i(t)$.
Let $D(t) = D_1(t)D_2(t) = \alpha_\kappa t^\kappa + \cdots +
\alpha_0$ ($\kappa \le k$), $c(t) = c_1(t)c_2(t)$.

Apply Lemma \ref{lemftgt} to $c(t)$ and $D_2(t)$. Thus $\exists
c'$ such that if $\exists t$ where $p$ divides both polynomials,
then $p|c'$. Let $\mathcal{P}_2$ be the prime divisors of $c'$ not
in $\mathcal{P}_0$ and let $\mathcal{P}_1$  be the prime divisors
of $\alpha_\kappa \cdot \mbox{Discriminant}(D(t))$ not in
$\mathcal{P}_0$. Define

\begin{eqnarray}
\mathcal{P} &=& \bigcup_{i=1}^2 \mathcal{P}_i, \ \ \ P = \prod_{p
\in \mathcal{P}} p. \end{eqnarray}

Note every prime in $\mathcal{P}$ is greater than $k$ and not in
$\mathcal{P}_0$.

As the product of primitive polynomials is primitive, $D(t)$ is
primitive. For any prime, either $D(t)$ mod $p$ is a constant not
divisible by $p$ or a non-constant polynomial of degree at most
$k$. In the second case, as there are at most $k$ roots to $D(t)
\equiv 0$ mod $p$, we find that given a $p > k$, $\exists t_p$
such that $D(t_p) \not\equiv 0$ mod $p$. By the Chinese Remainder
Theorem, $\exists t_0 \equiv t_p$ mod $p$ for all $p \in
\mathcal{P}$.

\subsubsection{Calculating the Conductor}

$\forall p \in \mathcal{P}$, $D(Pt + t_0) \equiv D(t_0) \not\equiv
0$ mod $p$. As $\mathcal{P}$ and $\mathcal{P}_0$ are disjoint,
this implies that $D(Pt + t_0)$ is minimal for all $p \in
\mathcal{P}$, as $\mathcal{P}_0$ contains the factors of $d$,$2$
and $3$. Moreover, $f_p(Pt+t_0) = 0$ for $p \in \mathcal{P}$.

By changing variables again, from $t$ to $P_0^mt + t_1$, we can
determine the powers of $p \in \mathcal{P}_0$ in the conductor.
Combining the two changes, we send $t$ to $\tau = P(P_0^m t + t_1)
+ t_0$.

Originally we had $\Delta(t) = d \Delta_1(t) \Delta_2(t)$. Now we
have $\Delta(\tau) = d \Delta_1(\tau) \Delta_2(\tau)$. It is
possible that $D_1(\tau) D_2(\tau)$ is no longer primitive;
however, if there is a common prime divisor $p$, $p$ divides
$\alpha_\kappa (P\cdot P_0^m)^\kappa$, implying $p \in
\mathcal{P}_0 \sqcup \mathcal{P}$.

We sieve to $D(\tau)$ square-free for $p \not\in \mathcal{P}_0
\sqcup \mathcal{P}$. As $\mathcal{P}_0 \sqcup \mathcal{P}$
contains all primes less than $k$, as well as the prime divisors
of $P_0$, $P$, $\alpha_\kappa$ and
$\mbox{Discriminant}(\Delta(t))$, we can perform the sieving. Note
the discriminants of $\Delta(t)$ and $\Delta(\tau)$ differ by a
power of $P\cdot P_0^m$. Thus, away from these primes, $D(\tau)
\equiv 0$ mod $p^2$ has at most $k < p^2$ roots, and we may sieve
to a positive percent of $t$. The sieving is unconditional if each
irreducible factor of $D(\tau)$ is of degree at most $3$.

$D(\tau)$ is divisible by fixed powers of primes in
$\mathcal{P}_0$ and never divisible by primes in $\mathcal{P}$.
Thus $\exists c_1$, $c_2$ with factors in $\mathcal{P}_0$ such
that $\widetilde{D}(\tau) = \frac{D_1(\tau)}{c_1}
\frac{D_2(\tau)}{c_2}$ is not divisible by any $p \in
\mathcal{P}_0 \sqcup \mathcal{P}$. We sieve to
$\widetilde{D}(\tau)$ square-free; for $p \not\in \mathcal{P}_0
\sqcup \mathcal{P}$, this is the same as $D(\tau)$ not divisible
by $p^2$.

We need to determine $f_p(\tau)$ for $p \in \mathcal{P}_0$, $p \in
\mathcal{P}$, and $p \not\in \mathcal{P}_0 \sqcup \mathcal{P}$.

By our previous arguments, if $m$ is sufficiently large,
$f_p(\tau) = f_p(Pt_1+t_0)$ for $p \in \mathcal{P}_0$.

If $p \in \mathcal{P}$ then $p \not\in \mathcal{P}_0$. Mod $p$,
$\Delta(\tau) = \Delta\Big( P(P_0t + t_1) + t_0\Big) \equiv
\Delta(t_0) \not\equiv 0$. Thus, for these $p$, $f_p(\tau) = 0$.

Assume $p \not\in \mathcal{P}_0 \sqcup \mathcal{P}$. The leading
term of $dD(\tau)$ is $d \alpha_\kappa (P\cdot P_0^m)^\kappa$. By
construction, $p$ does not divide the leading coefficient of
$\Delta(\tau)$, as $\mathcal{P}_0 \sqcup \mathcal{P}$ contains the
prime divisors of $d$, $\alpha_k$, $P$ and $P_0$. If we sieve to
$\widetilde{D}(\tau)$ square-free for $p \not\in \mathcal{P}_0
\sqcup \mathcal{P}$, then as the degree of $\Delta(\tau)$ is at
most $10$, the curve is minimal for such $p$. Thus, $f_p(\tau)$ is
$1$ if $p|D_2(\tau)$ and $2$ if $p|D_1(\tau)$.

Thus, we have shown

\begin{thm} All quantities as above, for $\widetilde{D}(\tau)$ square-free,
the conductors are

\begin{eqnarray}
C(\tau) = \prod_{p \in \mathcal{P}_0} p^{f_p} \cdot
\Bigg(\frac{|D_1(\tau)|}{c_1}\Bigg)^2 \frac{|D_2(\tau)|}{c_2}.
\end{eqnarray}

For sufficiently large $\tau$, $C(\tau)$ is a monotone increasing
polynomial (we may drop the absolute values), and a positive
percent of $\tau$ yield $\widetilde{D}(\tau)$ square-free.

\end{thm}

\subsection{Rational Surfaces II}

We consider what could go wrong in our proof if we are in case
$(2)$, where $3\mbox{deg}\ A(t) = 2\mbox{deg}\ B(t) = 12$ and
$\mbox{ord}_{t=0}t^{12} \Delta(t^{-1}) = 0$.

Thus, $\Delta(t)$ is a degree twelve polynomial, and we need to
worry about minimality issues. As before, we have

\begin{eqnarray}
\Delta(t) &=& -2^4\Big(2^2 A^3(t) + 3^3 B^2(t)\Big) = d
\Delta_1(t) \Delta_2(t) \nonumber\\ c_4(t) &=& c \gamma_1(t)
\gamma_2(t) \nonumber\\ \mathcal{P}_0 &=& \{p: p \leq \mbox{deg}\
\Delta(t) \} \ \cup \ \{p: p| cd\}, \ \ \ P_0 = \prod_{p \in
\mathcal{P}_0} p. \end{eqnarray}

There are three cases:

\bi

\item $\Delta(t)$ not divisible by a twelfth power;

\item $(\alpha t + \beta)^{12}|\Delta(t)$, $(\alpha t + \beta)
\notdiv c_4(t)$;

\item$(\alpha t + \beta)^{12}|\Delta(t)$, $(\alpha t + \beta) |
c_4(t)$.

\ei

These cases are handled in a similar fashion as before; see
\cite{Mil} for the calculations.

\subsection{Generalizations}

The previous arguments are applicable to any family where
$\mbox{deg}\ \Delta(t) \leq 12$ (which can include some
non-rational families). It is straightforward to generalize these
arguments for all families.

\subsection{Summary}

We summarize our sieving and conductor results:

\begin{thm}[Conductors and Cardinalities for
Families]\label{thmcondcard} For a one-parameter family with
$\mbox{deg}\ \Delta(t) \leq 12$, which includes all rational
families, by sieving to a positive percent subsequence we obtain a
family with conductors given by a monotone polynomial; further, by
Theorem \ref{thmconditionsnfbig}, after changing variables to
$\tau = P^m t + t_0$, a positive percent of $t \in [N,N]$ give
$D(\tau)$ square-free except for primes $p|P$, where the power of
such $p$ dividing $D(\tau)$ is independent of $t$. If all the
irreducible factors of $\Delta(t)$ are degree $3$ or less, the
sieving is unconditional; for degree $4$ and higher, the sieving
is a consequence of the ABC or Square-Free Sieve conjecture.
\end{thm}


\section{Sums of Test Functions at Primes}
\setcounter{equation}{0}

We list several standard sums of test functions over primes.
$\widehat{F}$, $\widehat{f}_i$ are even Schwartz functions with
compact support, $\varphi(m)$ is the Euler phi-function.

All statements below are straightforward applications of partial
summation and RH (or GRH for Dirichlet $L$-functions if $m \neq
1$) to handle the prime sums (see, for example, \cite{Mil});
weaker error terms are obtainable by the Prime Number Theorem.

\begin{lem}[Sum of $\widehat{F}$ over
primes]\label{lemprimeonesum}
\begin{eqnarray}
\frac{1}{\log N} \sum_{p \equiv b (m)} \frac{\log p}{p}
\widehat{F}\Big(a \flogpn\Big) = \frac{1}{2a\varphi(m)} F(0) +
\fologn. \end{eqnarray}
\end{lem}

Setting $m = 1$ and $a = 1,2$ yields

\begin{cor} \label{primeonesuma}$\frac{1}{\log N} \sum_p \frac{\log p}{p}
\widehat{F}\Big(\flogpn\Big) = \frac{1}{2} F(0) + \fologn.$
\end{cor}

\begin{cor} \label{primeonesumb}$\frac{1}{\log N} \sum_p \frac{\log p}{p}
\widehat{F}\Big(2\flogpn\Big) = \frac{1}{4} F(0) + \fologn.$
\end{cor}

\begin{lem}\label{primetwosuma}
\begin{eqnarray}
4 \sum_p \frac{\log^2 p}{\log^2 M} \frac{1}{p} \hfo
\hft\Big(\plogm\Big) &=& 2 \int_{-\infty}^{\infty} |u| \hfo(u)
\hft(u) du + O\Big(\frac{1}{\log M}\Big). \nonumber\\
\end{eqnarray}
\end{lem}

For $p \equiv b (m)$ we have

\begin{lem}\label{primetwosumb}
\begin{eqnarray}
4 \sum_{p \equiv b (m)} \frac{\log^2 p}{\log^2 M} \frac{1}{p} \hfo
\hft\Big(\plogm\Big) &=& \frac{2}{\varphi(m)}
\int_{-\infty}^{\infty} |u| \hfo(u) \hft(u) du +
O\Big(\frac{1}{\log M}\Big).\nonumber\\ \end{eqnarray}
\end{lem}

\begin{lem}\label{lemaepone}
Let $\mathcal{E}$ have rank $r$ over $\Qoft$ and assume Tate's
conjecture for $\mathcal{E}$ (known if $\mathcal{E}$ is a rational
surface). Then
\begin{eqnarray}
2\sum_p \plogx \oop \widehat{f}\Big(\plogx\Big)
\frac{-A_{1,\mathcal{F}}(p)}{p} = rf(0) \ + \ o(1). \end{eqnarray}
\end{lem}

Finally, we constantly encounter sums such as

\begin{eqnarray}
\sum_p \frac{\log p}{\log C(t)} \frac{1}{p^r}
\widehat{f}\Big(r\frac{\log p}{\log C(t)}\Big) a_t^r(p),
\end{eqnarray}

where $r \in \{1,2\}$ and $\log C(t)$ is $k \log N + o(\log N)$.

By Hasse, $a_t^r(p) \leq (2\sqrt{p})^r$. The contribution $S_l$
from $p \leq \log^l N$ is

\begin{eqnarray}
S_l & \ll & \frac{1}{\log N}\sum_{p \leq \log^l N} \frac{\log
p}{p^{r/2}}.
\end{eqnarray}

Clearly the larger contribution is from $r = 1$. By the Prime
Number Theorem, $\sum_{p \leq x}$ $\log p$ $\ll x$. By partial
summation, $\sum_{p \leq x} \frac{\log p}{\sqrt{p}}$ $\ll
\sqrt{x}$. Thus

\begin{eqnarray}
S_l \ll \frac{\sqrt{\log^l N} }{\log N}. \end{eqnarray}

We have shown

\begin{lem}[Removing Small Primes]\label{lemsmallprimes}
The sums over primes $p \leq \log^l N$ in the Explicit Formula
contribute $O(\log^{\frac{l}{2} - 1} N)$. For $l < 2$, this is
negligible.
\end{lem}


\section{Handling the Error Terms in the $2$-Level Density}
\setcounter{equation}{0}

Following Rudnick-Sarnak \cite{RS} and Rubinstein \cite{Ru}, we
handle the error terms in the $2$-level density, assuming we are
able to prove the $1$-level density theorem with error terms. By
the Explicit Formula (Equation \ref{thmef})

\begin{eqnarray}
\sum_{j_i} F_i\Big(\frac{\log N_E }{2\pi} \gamma_E^{(j_i)}\Big)
&=& \textbf{Good}_i + O\Big( (\log N_E)^{-\foh} \Big),
\end{eqnarray}

where \textbf{Good}$_i$ is the good part of the Explicit Formula,
involving $\widehat{F}(0)$, $F(0)$, and sums of $a_E(p)$ and
$a_E^2(p)$ for primes $p > \log N$.

Multiplying and summing over $i$ yields

\begin{eqnarray}
\oof \sumef \prod_{i=1}^{2} \Bigg[ \sum_{j_i} F_i\Big(\frac{\log
N_E }{2\pi} \gamma_E^{(j_i)}\Big) + O\Big( (\log N_E)^{-\foh}
\Big) \Bigg] &=& \oof \sumef  \prod_{i=1}^{2} \textbf{Good}_i.\nonumber\\
\end{eqnarray}

Multiplying out the LHS yields terms like

\begin{eqnarray}
O\Bigg[ \oof \sumef (\log N_E)^{-\frac{2-k}{2}}  \prod_{m=1}^{k}
\sum_{j_{m_i}} F_i\Big(\frac{\log N_E}{2\pi} \gamma_E^{(j_{m_i})}
\Big) \Bigg]. \end{eqnarray}

If each function $F_i$ were positive, we could insert absolute
values and move $\oof \sumef$ past the $\log^{-\frac{2-k}{2}} N_E$
factor. We assume our family has been sieved, so that the
conductors satisfy $\log N_E = c\log N + o(\log N)$.

There are three terms. If $k=0$ there is clearly no net
contribution. For $k = 1$ we have a $1$-level density, which is
finite by assumption. No error hits the $k=2$ piece (this is the
piece we want to calculate!). Only the $k=1$ piece is troublesome
for $F_i$ not positive.

If $F_i$ is not positive, we increase the above by replacing $F_i$
with a positive function $g_i$ such that $g_i$ is an even Schwartz
function whose Fourier Transform is supported in the same interval
as that of $F_i$ and $g_i(x) \geq |F_i(x)|$. As the $g_i$ satisfy
the necessary conditions, we may apply the $1$-Level Density
Theorem to the $g_i$'s, obtaining a bounded quantity. Hitting this
with $(\log N_E)^{-\foh}$, we see there is negligible
contribution.

For a construction of $g_i$, see Rubinstein \cite{Ru}, pages
$40-41$ or Rudnick-Sarnak \cite{RS}, pages $302-304$.

We have shown:

\begin{thm}[Handling the Error Terms]\label{thmerror}
If we are able to do the $1$-level density calculations, then we
may ignore the error terms in the $2$-level density.
\end{thm}

Note: the error need not be $O(\log^{-\foh} N)$; $o(1)$ also
works.\\ \ \\

\thanks{I wish to thank Harald Helfgott, Henryk Iwaniec, Nick
Katz, Wenzhi Luo, and Peter Sarnak for many enlightening
conversations.}

\end{document}